\documentclass[reqno]{amsart}
\usepackage{amsmath,amsfonts,amssymb,amsthm,graphicx,mathtools}
\usepackage[usenames,dvipsnames]{xcolor}
\usepackage[colorlinks=true,linkcolor=blue,urlcolor=blue]{hyperref}
\usepackage{marginnote}
\usepackage{overpic}
\newtheorem{thm}{Theorem}[section]
\newtheorem{lem}[thm]{Lemma}
\newtheorem{prop}[thm]{Proposition}
\newtheorem{cor}[thm]{Corollary}
\theoremstyle{definition}
\newtheorem{defn}[thm]{Definition}
\newtheorem{exmpl}[thm]{Example}
\newtheorem{rmk}[thm]{Remark}
\newtheorem{rmks}[thm]{Remarks}
\newtheorem{assu}[thm]{Assumption}
\numberwithin{equation}{section}

\newcommand{\N}{\mathbb{N}}
\newcommand{\R}{\mathbb{R}}
\newcommand{\Q}{\mathbb{Q}}

\newcommand{\B}{\mathcal{B}}

\newcommand{\K}{\mathcal{K}}
\newcommand{\M}{\mathcal{M}}
\newcommand{\T}{\mathcal{T}}
\newcommand{\U}{\mathcal{U}}
\newcommand{\Om}{\mathcal{O}}
\newcommand{\A}{\mathbb{A}}
\newcommand{\ep}{\varepsilon}
\newcommand{\wt}{\widetilde}

\newcommand{\PM}{\mathfrak{M_p}}
\newcommand{\nbd}{\nobreakdash}

\newcommand{\nti}{{n\to\infty}}

\DeclareMathOperator{\cls}{cls}

\DeclareMathOperator{\supp}{supp}

\begin{document}
\title[A new class of generalized ODEs]
{A new class of generalized ordinary differential equations with applications}
\address[I.P. Longo]{University of Exeter, Department of Mathematics and Statistics, 925
Laver Building, 23 N Park Road, EX4 4QE Exeter, United Kingdom.}
\email{i.longo@exeter.ac.uk}
\address[Sylvia Novo]{Departamento de Matem\'{a}tica Aplicada, E.I. Industriales, Universidad de
  Valladolid, Paseo Prado de la Magdalena 3-5
47011 Valladolid, Spain }
\email{sylvia.novo@uva.es}
\address[Rafael Obaya]{Departamento de Matem\'{a}tica Aplicada, E.I. Industriales, Universidad de
  Valladolid, Paseo Prado de la Magdalena 3-5
47011 Valladolid, Spain  and member of IMUVA, Instituto de Investigaci\'{o}n en
Matem\'{a}ticas, Universidad de Va\-lla\-dolid, Spain.}
\email{rafael.obaya@uva.es}
\address[Ana M. Sanz]{Departamento de Did\'{a}ctica de las Ciencias Experimentales, Sociales y de la Matem\'{a}tica,
Facultad de Educaci\'{o}n y Trabajo Social, Universidad de Valladolid, 47011 Valladolid, Spain,
and member of IMUVA, Instituto de Investigaci\'{o}n en  Mate\-m\'{a}\-ti\-cas, Universidad de
Valladolid, Spain.}
\email{anamaria.sanz@uva.es}
\thanks{Funding Declaration:~All authors were partly supported by MICIIN/FEDER projects PID2021-125446NB-I00 and PID2024-156691NB-I00. R. Obaya and A.M. Sanz were partly supported by JCyL/FEDER project CLU-2025-1-02-IMUVA. I.P. Longo was partly supported by Horizon Europe project ``Past to Future" (P2F), grant agreement number 101184070.}
\author[I.P.~Longo]{Iacopo P. Longo}
\author[S.~Novo]{Sylvia Novo}
\author[R.~Obaya]{Rafael Obaya}
\author[A. M.~Sanz]{Ana M. Sanz}
\date{}
\begin{abstract}
The paper presents a new approach to the study of precompact families of Carath\'{e}odory ODEs using techniques from nonautonomous dynamical systems. This problem falls outside the context of ordinary equations. We introduce the space of parametric b-measures, endowed with suitable topologies, in order to define a precise class of generalized ODEs. As an application, we study the dynamics of the fast variables in a slow–fast system of ODEs, where the fast motion is governed by a Carath\'{e}odory vector field with equicontinuous $m$-bounds and bounded $l$-bounds. Numerical simulations are provided to illustrate the scope and applicability of the conclusions.
\end{abstract}
\maketitle
\section{Introduction}\label{secintro}
In this paper, we introduce a precise class of generalized ordinary differential equations (ODEs) associated with parametric b-measures. The initial motivation arises from the observation that solutions of limiting equations in the Carath\'{e}odory ODE framework may lose the property of absolute continuity while remaining continuous and of bounded variation on compact intervals. Several approaches can be used to address these irregular limiting equations, and we propose a novel method. This class of equations offers an intermediate modeling framework between classical ODEs and stochastic equations. It is anticipated that they will find potential applications in fields such as financial markets, epidemiology, and climate science. Here, the application is devoted to the behaviour of the fast variables of a slow-fast system of ODEs, supported by some numerical examples in which the tracking phenomenon can be observed.

The study of Carath\'{e}odory ODEs is a classical question that has received intense attention since its inception, both for the interest in the theory itself and for the importance of its applications in science and engineering problems. Recently, Longo~\cite{tesis:Iacopo} and Longo et al.~\cite{paper:LNO,paper:LNO2} have introduced and investigated strong and weak topologies of integral type defined on appropriate spaces of Carath\'{e}odory functions. The convergence of a sequence $\{f_n\}_{n \geq 1}$, with respect to one of these topologies, requires the convergence of the integrals of the evaluation of the functions $f_n$ either pointwise over a dense and countable subset $D \subset \R^N$ (topologies $\T_D, \sigma_D$) or uniformly over sets of bounded continuous functions with a common modulus of continuity (topologies $\T_\Theta, \sigma_\Theta$). These works contain sufficient conditions on a set $E$ of Carath\'{e}odory functions implying that the mentioned topologies are topologies of continuity on $E$ for the corresponding Carath\'{e}odory ODEs.  Longo et al.~\cite{paper:LNuO,lno2} apply these conclusions to the study of mathematical models of critical transitions. Carath\'{e}odory ODEs also play a role in the theory of random ordinary differential equations (RODEs) if the noise process in a RODE is only measurable rather than continuous in time (see Han and Kloeden~\cite{book:HK}).

An important result in this theory states that, if $E$ is a set of Carath\'{e}odory functions with  equicontinuous $m$-bounds and bounded $l$-bounds and $D$ is a dense and countable subset of $\R^N$, then $\sigma_D$ is a topology of continuity on $E$ and the space $(E,\sigma_D)$ is a precompact metric space that, however, is not complete unless stronger conditions are satisfied (see Artstein~\cite{paper:ZA2} and~\cite{tesis:Iacopo}). In this paper, we provide a new construction for the compactification of $E$ by identifying the Carath\'{e}odory ODEs with a subset of the so-called {\it generalized ODEs given by parametric b-measures}  (a precise class of generalized ODEs governed by standard methods of ODEs), endowed with an extended $\sigma_D$ topology.  This construction simplifies the application of dynamical methods, based on the existence of ergodic measures, compact omega-limit sets and local, pullback or global attractors, in the study of precompact families of Carath\'{e}odory ODEs, and allows the extension to this context of previous results stated in the literature for more restrictive families of ODEs.

In fact, this particular class of generalized ODEs is part of the classical measure differential equations in the terms stated by Das and Sharma~\cite{paper:DSh},  Pandit and Deo~\cite{paper:PD} and Schmaedeke~\cite{paper:Schma}; and it also falls within the class of generalized  Kurzweil differential equations whose main properties can be found in Kurzweil~\cite{paper:kurz,book:kurz} and Schwabik~\cite{book:Schw} (the more recent reference Bonotto et al.~\cite{book:Bonotto} addresses a revised and expanded version of these questions with applications on Banach spaces). Some advantages of our construction are that these equations inherit their properties directly from the Carath\'{e}odory equations and that only Riemann integration and techniques of ODEs are needed.

Next, we describe the paper's organization and its main achievements. In Section~\ref{sec Motivation} we introduce the precompact families of Carath\'{e}odory ODEs and give a short motivation for the definition of the so-called parametric b-measures, which are fundamental in our approach. Inspired by the formula $d\nu_y = g(y,t)\,dt$ ($y\in\R$) in the scalar case,  a parametric b-measure on $\R^N$  is a parametric family $\nu \equiv \{\nu_y\}_{ y \in \R^N}$ of real Borel measures on bounded sets of $\R$ without purely discontinuous part and with $m$-bounds and $l$-bounds given by positive continuous measures.

The basic theory of parametric b-measures is established in Sections~\ref{sec Integ}--\ref{sec hull}. In Section~\ref{sec Integ}, the Riemann approach is used to define the integral of a parametric b-measure along a continuous curve over a compact interval, $\int_a^b d\nu_{y(s)}$.
Section~\ref{sectopo} contains extended versions of the weak topologies $\sigma_\Theta, \sigma_D$ mentioned above applied to the space $\PM$ of parametric b-measures on $\R^N$.
Some topological properties of a subset $E\subset \PM$, such as boundedness or equicontinuity of the $m$-bounds/$l$-bounds, are transferred to its closure in $\PM$, $\overline E:=\mathrm{cls}_{(\PM,\T)}(E)$ for  $\T\in\{\sigma_\Theta, \sigma_D\}$. Furthermore, if $E$ has bounded $l$-bounds, both topologies coincide. The main result asserts that if $E$ has equicontinuous $m$-bounds and  bounded $l$-bounds, then $\overline E$ is compact. Thus, for a precompact family of Carath\'{e}odory functions, the  closure of its embedding into the space of parametric b-measures defines its compactification.

In Section~\ref{sec hull} we address the continuity of the time translation map  $\R \times \overline E \to \PM$, $(t,\nu) \mapsto \nu {\cdot} t$. If $\nu$ is a parametric b-measure on $\R^N$ with equicontinuous $m$-bounds,  then the hull of $\nu$, $\mathrm{Hull}_{(\PM,\T)}(\nu):=\mathrm{cls}_{(\PM,\T)}\{\nu{\cdot}t\mid t\in\R\}$ is invariant under translation and $(t,\nu) \mapsto \nu {\cdot} t$ defines a continuous flow on it, for $\T\in\{\sigma_\Theta, \sigma_D\}$.

The second part of the paper is devoted to the study of generalized ODEs $y'= \bar\nu_{y(t)}$ given by parametric b-measures and their application to Carath\'{e}odory ODEs. These equations are introduced and investigated in Section~\ref{sec generalized ODEs}. Solutions are continuous functions that satisfy the integral equation $y(t)=y(t_0)+\int_{t_0}^td\bar\nu_{y(s)}$; hence, they are continuous functions of bounded variation on compact intervals, but in general they are not absolutely continuous. Standard techniques of ODEs allow us to state the basic theory on local existence, uniqueness, and extension of solutions. Also in this context   $\sigma_D$ is a topology of continuity for the solutions under the assumptions of equicontinuous $m$-bounds and bounded $l$-bounds. As a consequence, if $\bar\nu_0 \in \PM^{\!\!\!N}$ has the former properties,  a local continuous skew-product flow on $\mathrm{Hull}_{(\PM^{\!\!\!N},\bar\sigma_D)}(\bar\nu_0) \times \R^N$ is defined through translation in the hull and the solutions of the generalized ODEs over the hull in the fiber.
All these results are applied in Section~\ref{subsec appl Carath}, where a Carath\'{e}odory ODE is identified with a generalized ODE given by an absolutely continuous parametric b-measure, preserving the solutions. 
Some numerical examples are provided in Section~\ref{subsec:numerics}, to show how solutions of the type of differential equation described in this work might look like.

Finally, Section~\ref{sec slowfast} illustrates an application of the previous theory in a context of slow-fast systems of ODEs, where the fast motion is governed by a Carath\'{e}odory map with equicontinuous $m$-bounds and bounded $l$-bounds. By identifying the equations with generalized ODEs given by parametric b-measures, we can build the hull, which is compact, and the skew-product semiflow induced by the solutions. Then, techniques of nonautonomous dynamical systems allow us to extend previous results by Artstein~\cite{paper:ZA99} and Longo et al.~\cite{paper:LOStracking}, now under rather relaxed conditions. We complete the section with an analytical and numerical study of a scalar ordinary differential equation featuring an explicit dependence on time with two distinct time scales---the equation can be rewritten as a slow-fast system of the type we study in the rest of the section. We showcase that solutions of the slow-fast equation track (on suitable intervals) a family of solutions of nonautonomous equations with sole dependence on the fast time-scale, obtained from the original slow-fast system fixing the slow changing coefficients as parameters.

For the sake of completeness, we finish the introduction with some preliminaries of measure theory. As usual,  $\R^N$  is the $N$-dimensional Euclidean space with norm $|\cdot|$  and $B_r$ is the closed ball of $\R^N$ centered at the origin and with radius $r$. $C_c(\R)$ will denote the space of real continuous functions on the real line with compact support endowed with the norm $\|\cdot\|_\infty$. We denote by  $\mathcal{A}$ the $\sigma$-algebra of Borel sets in $\R$ and by $\M^+$  the set of  locally finite regular Borel (positive) measures on $\R$ endowed with the vague topology, i.e., a sequence  $\{\mu_n\}_{n\geq 1}$ of measures in $\M^+$ \emph{vaguely converges to $\mu\in\M^+$}, and we write it as $\mu_n\xrightarrow[]{\widetilde{\sigma}}\mu$, if and only if
\begin{equation*}
 \lim_{\nti}\int_\R \phi(s)\, d\mu_n(s)=\int_\R\phi(s)\, d\mu(s)\qquad \text{for each } \phi\in C_c(\R)\,.
\end{equation*}
As shown in Kallenberg~\cite[Theorem~15.7.7, p.170]{book:Kall}, $(\M^+,\widetilde{\sigma})$ is a Polish space, i.e., it is separable and completely metrizable.

From the decomposition theorem for regular Borel measures on $\R$ (see Hewitt and Stromberg~\cite[Theorem~19.61, p.337]{book:HS}), the set $\M^+$ decomposes as $\M^+=\M^+_{ac}\oplus\M^+_{sc}\oplus\M^+_{pd}$, the classes of absolutely continuous and singular continuous measures with respect to the Lebesgue measure $l$, and purely discontinuous measures,  respectively. We will denote by $\M^+_c$ the subset of measures without purely discontinuous part, that is,  $\M^+_c= \M^+_{ac}\oplus\M^+_{sc}$.

\section{Precompact families of Carath\'{e}odory ODEs and the motivation to introduce parametric b-measures}\label{sec Motivation}
In this paper, we consider Carath\'{e}odory ODEs $y'=g(y,t)$ with equicontinuous $m$-bounds and bounded $l$-bounds. More precisely,  $g\colon\R^N\times\R\to \R^N$ is (Lebesgue) measurable in $t$ for each $y\in\R^N$, it is continuous in $y$ for almost every $t\in\R$,  and for every compact set $K\subset \R^N$ there exist two nonnegative and locally integrable maps defined on $\R$, $\wt m_K, \wt l_K \in L^1_{loc}$ (the space of real Lebesgue locally integrable maps on $\R$), respectively called an $m$-bound and an $l$-bound of $g$ on $K$, such that:
\begin{itemize}
    \item[(i)] for a.e.~$t\in\R$,  $|g(y,t)| \leq  \wt m_K (t)$ for all $y\in K$;
    \item[(ii)] for a.e.~$t\in\R$,  $|g(y_1,t)-g(y_2,t)| \leq  \wt l_K (t)\,|y_1-y_2|$ for all $y_1,y_2\in K$.
\end{itemize}
The hypotheses of equicontinuous $m$-bounds and bounded $l$-bounds read as follows (see Artstein~\cite{paper:ZA2}): for every compact set $K\subset \R^N$,
\begin{itemize}
\item[(i)$'$] $\int_0^s \wt m_K (t)\,dt$ is uniformly continuous in $s$, that is, given $\ep>0$ there is $\delta=\delta_K(\ep)>0$ such that for all $s\in\R$, $\big|\!\int_s^{s+h} \wt  m_K (t)\,dt\big|<\ep$ if $|h|<\delta$;
\item[(ii)$'$] there exists a $c_K>0$ such that $\int_s^{s+1} \wt l_K (t)\,dt \leq c_K$ for all $s\in\R$.
\end{itemize}
As explained in~\cite{paper:ZA2} (see also~\cite{tesis:Iacopo}), under these conditions a challenging phenomenon appears: the limiting equations (i.e., the limits of a family of translated equations under a certain topology) may not be ODEs, that is, the family of translated equations is a precompact family but not necessarily complete. In order to complete it, a class of Kurzweil equations was considered in~\cite{paper:ZA2}, with the drawback that the properties of equicontinuous $m$-bounds and bounded $l$-bounds are lost for some of these equations, and the corresponding solutions can lose properties such as being of bounded variation in compact sets. Trying to maintain the properties of both the initial equation and its solutions in the compactificacion process, in Section~\ref{sec generalized ODEs} we propose a class of generalized ODEs, built from a measure-theoretic approach.

For a better understanding of our proposal, let us first consider a scalar linear equation $y'=a(t)\,y$, where $a(t)$ is a nonnegative locally integrable map defined on $\R$ such that, given $\ep>0$, there is $\delta>0$ such that $\left|\int_s^{s+h} a(t)\,dt\right|<\ep$ if  $|h|<\delta$ for all $s\in\R$. In this case, if with each translated map $a(t+s)$ ($s\in\R$) we associate the measure  $d(\mu{\cdot}s)=a(t+s)\,dt \in \mathcal{M}^+_{ac}$ ($s\in\R$), it holds that the set $M=\{\mu{\cdot}s\mid s\in \R\}$ is relatively compact in $( \mathcal{M}^+, \wt \sigma)$ and its closure keeps within the set of continuous measures (see~\cite[Theorem~2.23]{tesis:Iacopo}):
\begin{equation}\label{eq:cls continuas}
\cls_{(\mathcal{M}^+, \wt \sigma)}(M)\subset \mathcal{M}^+_{ac} \oplus \mathcal{M}^+_{sc}\,.
\end{equation}
Then, when we address the scalar nonlinear case $y'=g(y,t)$, for each fixed (parameter) $y\in\R$ one can consider the object $d\nu_y=g(y,t)\,dt$. (We advance that in Section~\ref{subsec appl Carath} we will rewrite the equation as $y'=\nu_{y(t)}$). Note that $d\nu_y=g(y,t)\,dt$ no longer defines a measure if $g$ is real and only locally integrable, but it defines a real measure when restricted to compact sets of $\R$. Also, \eqref{eq:cls continuas} suggests that we can hope to get a compact closure within a class of continuous objects. With this motivation, we define the concepts of (continuous real) \emph{b-measure} and \emph{parametric b-measure}, where the continuity condition is essential for our purposes.

Let $\mathcal{B}$ denote the ring of bounded Borel sets in $\R$ and let $\M_c(\B)$ denote the set of maps $\mu\colon \B\to\R$ that satisfy:
\begin{itemize}
  \item[(i)] $\mu(\emptyset)=0$;
  \item[(ii)]  $\mu(\cup_{n=1}^\infty A_n) = \sum_{n=1}^\infty \mu(A_n)$ for every countable family $\{A_n\}_{n\geq 1}\subset \B$ of pairwise disjoint sets such that $\cup_{n=1}^\infty A_n\in\B$;
  \item[(iii)] $\mu(\{t\})=0$ for every $t\in\R$.
\end{itemize}
\begin{defn}
The elements of $\M_c(\B)$ are called (continuous real) {\it b-measures}.
\end{defn}

\begin{rmks}\label{rm:b-measure}
(1) We write ``b-measure" to indicate that $\mu$ is defined on the ring of bounded Borel sets, where the countably additive property (ii) holds. Note that property (ii) implicitly requires that the series $\sum_{n=1}^\infty \mu(A_n)$ is unconditionally convergent, or equivalently, by Riemann's theorem,  absolutely convergent.
(2)~Property (iii) means that $\mu$ is continuous as a measure.  (3) When a b-measure $\mu$ is restricted to a compact set $K$ of $\R$, it defines a continuous real measure on~$K$; see~\cite[Section~19]{book:HS} for the theory of abstract integration with respect to real measures. (4) The set $\M_c(\B)$ has a vector structure: if $\mu_1,\mu_2\in \M_c(\B)$, the sum $(\mu_1+\mu_2)(A):=\mu_1(A)+\mu_2(A)$, $A\in \B$ defines a b-measure. Also, if $\alpha\in\R$, then $(\alpha\mu_1)(A):=\alpha\mu_1(A)$, $A\in \B$ defines a b-measure.
\end{rmks}

\begin{exmpl}
The set of b-measures includes two important cases. (i) The restriction to $\B$ of a continuous regular Borel signed measure is included in $\M_c(\B)$. (ii) If $f\in L^1_{loc}$, then the map  $\mu(A):=\int_A f(t)\,dt$, $A\in \B$ defines a b-measure.
\end{exmpl}

For a regular Borel signed measure $\mu$, one can define its {\it total variation} $|\mu|$ by  $|\mu|:=\mu^++\mu^-$ for the positive measures $\mu^+$ and $\mu^-$ given in its Jordan decomposition. In fact, as proved in~\cite[Theorem~19.10]{book:HS},
\begin{equation}\label{eq:total var}
|\mu|(A)=   \sup\bigg\{ \sum_{i=1}^n |\mu(A_i)|\, \bigg| \, \{A_1,\ldots,A_n\} \;\hbox{is a Borel partition of}\; A \bigg\}\,.
\end{equation}
Note that we can use this formula to define a positive b-measure $|\mu|$ from a real b-measure $\mu\in \M_c(\B)$, by taking partitions $\{A_1,\ldots,A_n\} \subset \B$.

We define what we mean by a parametric b-measure with respect to a non-decreasing family of moduli of continuity $\{\omega_j\colon\R^+\to\R^+\mid j\geq 1\}$ (that is, $\omega_j$ is continuous and non-decreasing, $\omega_j(0)=0$ and $\omega_j\leq \omega_{j+1}$ for all $j\geq 1$).
\begin{defn}\label{defi:parameasure}
Given a non-decreasing family of moduli of continuity $\{\omega_j\}_{j\geq 1}$, a \emph{parametric b-measure on $\R^N$} is a  map $\nu\colon \R^N\to \M_c(\B),\;y\mapsto \nu_y$
such that for each integer $j\geq 1$ there are positive continuous measures $m_j, l_j\in \M_c^+$ satisfying:
\begin{itemize}
\item[(i)] $|\nu_y|(A)\leq m_j(A)$ for each subset $A\in \B$ and each  $y\in\R^N$ with $|y|\leq j$;
\item[(ii)] $|\nu_{y_1}-\nu_{y_2}|(A)\leq l_j(A) \,\omega_j(|y_1-y_2|)$ for each subset $A\in\B$ and $y_1, y_2\in\R^N$ with $|y_1|, |y_2|\leq j$.
\end{itemize}
The measures $m_j$ and $l_j$ of $\M_c^+$ are called \emph{$m$-bounds} and \emph{$l$-bounds} of $\nu$ on the closed ball $B_j$, respectively. The set of parametric b-measures on $\R^N$ with respect to the family of moduli $\{\omega_j\}_{j\geq 1}$ is denoted by $\PM$. An \emph{$M$-dim  parametric b-measure on $\R^N$} is a  map $\bar\nu\colon \R^N\to \M_c(\B)^M,\;y\mapsto (\nu^1_y,\ldots,\nu^M_y)$ such that each component $\nu^i$, $1\leq i\leq M$ is a parametric b-measure on $\R^N$.
\end{defn}
Note that $\{\nu_y\}_{y\in\R^N}$ is a parametric family of b-measures and $y$ is a vector-valued parameter. When we look at the map $\nu\colon\R^N\to \M_c(\B)$, (ii) is a continuity property of $\nu$ on $B_j$ given through the modulus of continuity $\omega_j$. The name \emph{$l$-bound} comes from the standard Lipschitz-type condition in the literature of Carath\'{e}odory differential equations (see, e.g.,~\cite{paper:ZA2,paper:LNO,paper:LNO2}). Here we keep that name, although we admit other moduli of uniform continuity rather than the identity map.

It is easy to check that the space $\PM$  of parametric b-measures on $\R^N$ with respect to a family of moduli of continuity $\{\omega_j\}_{j\geq 1}$ has a vector structure with standard definitions.
Let us give a simple example for later use.
\begin{prop}\label{prop:phinu}
Given a b-measure $\nu_0\in \M_c(\B)$ that has an m-bound $m_0\in\M^+_c$ over $\B$ and a map $\phi\colon \R^N\to[0,1]$ of class $C^1$ with compact support, the  map $\,\phi\,\nu_{0}\colon \R^N\mapsto  \M_c(\B) $, $y\to \phi(y)\nu_{0}$
defines a parametric b-measure on $\R^N$ with respect to the trivial family of moduli given by the identity map.
\end{prop}
\begin{proof}
Clearly $\phi\,\nu_{0}$ is well-defined.
Note that $|\phi(y)\nu_{0}|(A) = \phi(y)\, |\nu_{0}|(A) \leq m_0(A)$ for each subset $A\in \B$ and every $y\in\R^N$. Moreover, for all $y_1, y_2\in\R^N$, by the Mean Value theorem there exists a $c_0>0$ such that
\[
|\phi(y_1)\nu_{0}-\phi(y_2)\nu_{0}|(A) = |\phi(y_1)-\phi(y_2)|\,|\nu_{0}|(A) \leq c_0\,m_0(A)\,|y_1-y_2|\,.
\]
Therefore, (i) and (ii) of Definition~\ref{defi:parameasure} hold with $m$-bounds and $l$-bounds, and modulus $\omega_j = {\rm id}$, all independent of $j\geq 1$. The proof is finished.
\end{proof}

\section{Integration of parametric b-measures along curves}\label{sec Integ}
Since we will eventually solve equations given by parametric b-measures, we need an integration theory, and the Riemann approach is enough for our purposes. To introduce the integration of real maps with respect to parametric b-measures evaluated along curves $y\colon[a,b]\to \R^N$, we give the following definitions.
\begin{defn}\label{defi:division}
	  Let $\delta >0$.  A \emph{tagged} $\delta$\nbd-\emph{partition} of the real interval $[a,b]$ is a set
	  \[
	  \Delta^\delta\coloneqq\{(\tau_i, [t_i,t_{i+1}])\mid i=1,\ldots, k-1\}
	  \]
	   satisfying  $[a,b]=\bigcup_{i=1}^{k-1} [t_i,t_{i+1}]$, $a=t_1<t_2<\cdots<t_k=b$, $\tau_i\in [t_i,t_{i+1}]$,  and $t_{i+1}-t_i<\delta$ for each $i=1,\ldots,k-1$.
 Given a map $h\colon[a,b]\to \R$, a parametric b-measure $\nu$ on $\R^N$ and a curve $y\colon [a,b]\to \R^N$, we will denote by
	   \[
       S_\nu^{y(\cdot)}\big(h,\Delta^\delta\big)\coloneqq\sum_{i=1}^{k-1} \int_{t_i}^{t_{i+1}} h(\tau_i)\, d\nu_{y(\tau_i)}=\sum_{i=1}^{k-1} h(\tau_i)\, \nu_{y(\tau_i)}[t_i,t_{i+1}]
       \]
the (Riemann) sum corresponding to the tagged $\delta$-partition $\Delta^\delta$.
\end{defn}
\begin{defn}
Let $\nu$ be a parametric b-measure on $\R^N$ and consider a map $y\colon [a,b]\to \R^N$. A bounded map $h\colon [a,b]\to \R$ is  \emph{integrable with respect to $\nu$ evaluated along the curve $y\colon [a,b]\to \R^N$} if there exists an $L\in\R$ such that for each $\ep>0$ there is a $\delta(\ep)>0$ such that for every tagged $\delta(\ep)$\nbd-partition $\Delta^{\delta(\ep)}=\{(\tau_i, [t_i,t_{i+1}])\mid i=1,\ldots, k-1\}$ of $[a,b]$,
\begin{equation*}
\Big| L-S_\nu^{y(\cdot)}\big(h,\Delta^{\delta(\ep)}\big)\Big|=\bigg|L - \sum_{i=1}^{k-1} \int_{t_i}^{t_{i+1}} h(\tau_i)\,d\nu_{y(\tau_i)}\bigg|<\ep\,.
\end{equation*}
In this case, we write $\int_a^b h(s)\, d\nu_{y(s)}=L$. When the map $h$ is identically equal to $1$, we will simply say that the parametric b-measure $\nu$ is \emph{integrable along the curve} $y\colon [a,b]\to \R^N$ and write $\int_a^b  d\nu_{y(s)}$ for the integral and
\[
S_\nu^{y(\cdot)}\big(\Delta^\delta\big)= \sum_{i=1}^{k-1} \int_{t_i}^{t_{i+1}} d\nu_{y(\tau_i)} = \sum_{i=1}^{k-1} \nu_{y(\tau_i)}[t_i,t_{i+1}]
\]
for the Riemann sums. If $a<b$, then $\int_b^a  d\nu_{y(s)} := -\int_a^b  d\nu_{y(s)}$.
\end{defn}
As usual, as stated in the next result whose proof is omitted, the integrability can be checked in terms of the Cauchy property for the Riemann sums.
\begin{prop}\label{prop:integrability}
 The following statements are equivalent:
\begin{itemize}
\item[(i)] the map $h\colon [a,b]\to \R$ is integrable with respect to the parametric b-measure $\nu$ on $\R^N$ evaluated along the curve $y\colon [a,b]\to \R^N$;
\item[(ii)] given $\ep>0$ there is $\delta(\ep)>0$ such that $\big|S_\nu^{y(\cdot)}(h,\Delta_1^{\delta(\ep)})- S_\nu^{y(\cdot)}(h,\Delta_2^{\delta(\ep)})\big| <\ep$
for all tagged $\delta(\ep)$-partitions $\Delta_1^{\delta(\ep)}$ and $\,\Delta_2^{\delta(\ep)}$ of $[a,b]$.
\end{itemize}
\end{prop}
We do not include an exhaustive list with the standard properties of the integral defined in this way, but we remark that it is important that the b-measures $\nu_y$ are continuous, that is, they satisfy $\nu_y(\{t\})=0$ for all $t\in\R$ and $y\in\R^N$.

\begin{rmk}
To prove that every continuous map $h(t)$ is integrable with respect to a parametric b-measure $\nu$ evaluated along a continuous curve $y(t)$, it will suffice to study the integrability of the map $h_0(t)\equiv 1$, $t\in [a,b]$ along continuous curves. The reason is that one can build the parametric b-measure on $\R^{N+1}$, $\tilde \nu\colon\R^N\times \R\to \M_c(\B)$, $(y,y_{N+1})\mapsto y_{N+1}\,\nu_y$ (with respect to different moduli of continuity) and the continuous map $\tilde y\colon[a,b]\to \R^N\times \R$, $t\mapsto (y(t),h(t))$   and, given a tagged $\delta$-partition $\Delta^\delta$ of $[a,b]$, it holds that $S_\nu^{y(\cdot)}\big(h,\Delta^\delta\big) = S_{\tilde\nu}^{\tilde y(\cdot)}\big(\Delta^\delta\big)$.
\end{rmk}

To prove that each parametric b-measure is integrable along a continuous curve, we need the following technical lemma.
\begin{lem}\label{lem:tau}
 Given a parametric b-measure  $\nu\in\PM$ on $\R^N$ with respect to the family of moduli $\{\omega_j\}_{j\geq 1}$, a continuous map $y\colon [a,b]\to \R^N$  and an $\ep>0$,  there is a $\delta(\ep)>0$ such that if $|\tau_i-\wt\tau_i|<\delta(\ep)$ $(\tau_i,\wt\tau_i\in [a,b], 1\leq i\leq k-1)$, then
\[
\bigg|\,\sum_{i=1}^{k-1} \int_{t_i}^{t_{i+1}}\left( d\nu_{y(\tau_i)} - d\nu_{y(\wt\tau_i)}\right)\bigg|< \ep\,
\]
for every partition of  $[a,b]=\bigcup_{i=1}^{k-1} I_i$ with $I_i=[t_i,t_{i+1}]$ and $a=t_1<\cdots<t_k=b$.
\end{lem}
\begin{proof} Since $y$ is continuous in $[a,b]$, there is a $j\in\N$ such that $y([a,b])\subset B_j$. Let $l_j$  be an $l$-bound of $\nu$ on $B_j$, as in Definition~\ref{defi:parameasure}~(ii). Notice that $l_j([a,b])=\sum_{i=1}^{k-1} l_j(I_i)$ because $l_j\in \M_c^+$ has no purely discontinuous part. Denote $l_j ([a,b])= M $ and assume first that $M>0$. From the uniform continuity of $y$ and the continuity of the modulus $\omega_j$ at $0$, given $\ep>0$ there is a $\delta(\ep)>0$ such that $\omega_j(|y(\tau_i)-y(\wt\tau_i)|)< \ep/M$ provided that $|\tau_i-\wt\tau_i|<\delta(\ep)$. As a consequence,
\begin{align*}
&\bigg|\,\sum_{i=1}^{k-1}  \int_{t_i}^{t_{i+1}} \left( d\nu_{y(\tau_i)} - d\nu_{y(\wt\tau_i)}\right)\bigg|
 \leq \sum_{i=1}^{k-1} |\nu_{y(\tau_i)}- \nu_{y(\wt\tau_i)}|(I_i)\\ &\hspace{1cm}\leq \sum_{i=1}^{k-1} l_j(I_i)\, \omega_j(|y(\tau_i) - y(\wt\tau_i)|)
  < \sum_{i=1}^{k-1} l_j(I_i) \,\frac{\ep}{M}= l_j([a,b])\, \frac{\ep}{M}=\ep\,,
\end{align*}
as claimed. Finally, if $M=0$, the value of the sum is $0$ by the previous bounds.
\end{proof}
\begin{prop}\label{prop:defIntegral}
Every parametric b-measure $\nu$ on $\R^N$ is integrable along every continuous map $y\colon [a,b]\to \R^N$.
\end{prop}
\begin{proof} It is enough to check that (ii) in Proposition~\ref{prop:integrability} holds. Given $\ep>0$ we take  $\delta=\delta(\ep/2)>0$ provided in  Lemma~\ref{lem:tau} and two tagged $\delta$-partitions $\Delta_1^\delta$ and $\,\Delta_2^\delta$ of $[a,b]$.  By ordering all the points in these two partitions and adding some extra endpoints and tags, when necessary, we can build a new tagged $\delta$-partition $\Delta_3^\delta$ that is finer than both $\Delta_1^\delta$ and $\Delta_2^\delta$. Then,
\begin{equation*}
 \left|S_\nu^{y(\cdot)}(\Delta_1^\delta)- S_\nu^{y(\cdot)}(\Delta_2^\delta)\right|\leq \left|S_\nu^{y(\cdot)}(\Delta_1^\delta)- S_\nu^{y(\cdot)}(\Delta_3^\delta)\right| + \left|S_\nu^{y(\cdot)}(\Delta_3^\delta)- S_\nu^{y(\cdot)}(\Delta_2^\delta)\right|.
\end{equation*}
From the way of constructing $\Delta_3^\delta$ and Lemma~\ref{lem:tau} we can conclude that $\big|S_\nu^{y(\cdot)}(\Delta_1^\delta)- S_\nu^{y(\cdot)}(\Delta_2^\delta)\big|<\ep$ and the proof is finished.
\end{proof}
It is immediate to check the following result regarding the integration of a parametric b-measure of the form  $\phi\,\nu_{0}$, and the proof is omitted.
\begin{prop}\label{prop:Phinuintegral}
Let $\nu^*=\phi\,\nu_{0}$ be a parametric b-measure in the class considered in {\rm Proposition~\ref{prop:phinu}} and let $y\colon [a,b]\to \R^N$ be a continuous map. Then,
\[ \int_a^b d \nu^*_{y(s)}=\int_a^b \phi(y(s))\, d\nu_{0}\,. \]
\end{prop}
Note that $\int_a^b \phi(y(s))\, d\nu_{0}$ is the integral on $[a,b]$ of $\phi(y(s))$ with respect to the real measure ${\nu_0}|_{[a,b]}$, the restriction of $\nu_0$ to $[a,b]$ (see Remark~\ref{rm:b-measure}~(3)).

We finish this section with a useful result and an immediate corollary.
\begin{lem}\label{CotasInt}
Let $\nu\in\PM$ be a parametric b-measure on $\R^N$ with respect to the family of moduli $\{\omega_j\}_{j\geq 1}$, with $m$-bounds $\{m_j\}_{j\geq 1}$ and $l$-bounds $\{l_j\}_{j\geq 1}$, and consider two continuous maps $x,\,y\colon [a,b]\to \R^N$.  If $x([a,b]),y([a,b]) \subset B_j$, then
\[
\bigg|\int_a^b d\nu_{y(s)}\bigg|\leq m_j[a,b] \;\text{ and }\;  \bigg| \int_a^b [d\nu_{x(s)} - d\nu_{y(s)}]\bigg|\leq l_j[a,b]\,\omega_j(\|{x-y}\|_{[a,b]})\,,
\]
where $\|\,{\cdot}\,\|_{[a,b]}$ denotes the sup-norm on the interval $[a,b]$.
\end{lem}
\begin{proof}
Since the integral is the limit of the Riemann sums, for the first inequality, it suffices to check that for every tagged $\delta$-partition $\Delta^\delta=\{(\tau_i, [t_i,t_{i+1}])\mid i=1,\ldots, k-1\}$ of $[a,b]$,
$ |\sum_{i=1}^{k-1} \nu_{y(\tau_i)}(I_i)|
\leq \sum_{i=1}^{k-1} |\nu_{y(\tau_i)}|(I_i) \leq \sum_{i=1}^{k-1} m_j(I_i)= m_j[a,b]$,  where we have written $I_i=[t_i,t_{i+1}]$.  Analogously, from
\[ \bigg|\sum_{i=1}^{k-1}(\nu_{x(\tau_i)}-\nu_{y(\tau_i)})(I_i)\bigg| \leq \sum_{i=1}^{k-1}l_j(I_i)\,\omega_j(|x(\tau_i)-y(\tau_i)|)\leq l_j[a,b]\,\omega_j(\|{x-y}\|_{[a,b]})\,, \]
the second inequality is easily deduced. Notice that $m_j$, $l_j\in \M_c^+$ has been used.
\end{proof}

\begin{lem}\label{lem:xntox0}
Let $\nu\in \PM$ and let $\{y_n(\cdot)\}_{n\geq 1}\subset C([a,b],\R^N)$ be a sequence  converging uniformly  to $y_0(\cdot)$ in $[a,b]$.
Then,
 \[
 \lim_{n\to\infty}\int_a^b d\nu_{y_n(s)}=\int_a^b d\nu_{y_0(s)}\,.
 \]
\end{lem}
\begin{proof}
The proof follows from Lemma \ref{CotasInt}, once noted that by the uniform convergence there is a $j\geq 1$ such that $y_n([a,b]) \subset B_j$ for every $n\geq 0$.
\end{proof}

\section{Topologies in the space of parametric b-measures}\label{sectopo}
In this section, we introduce some topologies in the space  $\PM$ of parametric b-measures on $\R^N$  with respect to a fixed family of moduli of continuity  $\{\omega_j\}_{j\geq 1}$ (see Definition~\ref{defi:parameasure}). Then, we study how some topological properties of the $m$-bounds and $l$-bounds of a subset $E\subset \PM$ are transferred to its closure. The main result is the compactness of $E\subset \PM$ with equicontinuous $m$-bounds and bounded $l$-bounds.

In order to introduce a useful topology in $\PM$ based on suitable sets of moduli of continuity, as in~Longo et al.~\cite{paper:LNO, paper:LNO2} for spaces of Carath\'{e}odory functions,  we recall the definition.  For $I=[q_1,q_2]$, $ q_1,q_2\in\Q$ it is understood that $q_1< q_2$.
~\begin{defn}\label{def:modcont}
We call \emph{a suitable set of moduli of continuity} any countable  set of non-decreasing continuous functions
\begin{equation*}
\Theta=\left\{\theta^I_j \in C(\R^+, \R^+)\mid j\in\N, \ I=[q_1,q_2], \ q_1,q_2\in\Q\right\}
\end{equation*}
such that $\theta^I_j(0)=0$ for every $\theta^I_j\in\Theta$, and  with the relation of partial order given~by
\begin{equation*}\label{def:modCont}
\theta^{I_1}_{j_1}\le\theta^{I_2}_{j_2}\quad \text{whenever } I_1\subseteq I_2 \text{ and } j_1\le j_2 \,.
\end{equation*}
\end{defn}
We endow the space $\PM$ of parametric b-measures with the following topologies.
\begin{defn}[Topology $\sigma_\Theta$]\label{def:Tsigmatheta}
Let $\Theta$ be a suitable set of moduli of continuity. We call $\sigma_{\Theta}$ the topology on $\PM$ generated by the family of seminorms \vspace{-.025cm}
\begin{equation*}
p_{I, j}(\nu)=\sup_{y(\cdot)\in\K_j^I}\left|\,\int_I d\nu_{y(s)}\right| ,\quad \nu\in\PM\,,
\vspace{-.025cm}
\end{equation*}
for each $I=[q_1,q_2]$, $q_1,q_2\in\Q$ and $j\in\N$, where $\K_j^I$ denotes the set of functions $y(\cdot)$ in $C(I,B_j)$ which have $\theta^I_j$ as modulus of continuity, that is,
\begin{equation}\label{eq:KjI}
 \K_j^I:=\{y(\cdot)\in C(I,B_j)\mid |y(t_1)-y(t_2)|\leq \theta^I_j(|t_1-t_2|)\; \hbox{for all } t_1,t_2\in I\}\,.
\end{equation}
$\left(\PM,\sigma_{\Theta}\right)$ is a locally convex metric space.
\end{defn}
\begin{defn} [Topology $\sigma_D$]\label{def:TsigmaD}
Let $D$ be a dense and countable subset of $\R^N$. We call $\sigma_{D}$ the topology on $\PM$ generated by the family of seminorms
\begin{equation*}
p_{I, y}(\nu)=\left|\,\int_I d\nu_y\,\right|,\quad \nu\in\PM\,,
\end{equation*}
for  $y\in D$, $I=[q_1,q_2]$, $ q_1,q_2\in\Q$. $\left(\PM,\sigma_{D}\right)$ is a  locally convex metric space.
\end{defn}
\begin{rmk}\label{rmk:orden top}
Notice that $\sigma_D\leq \sigma_\Theta$.
\end{rmk}

If a function $y(\cdot)$ belongs to the set $\K_j^I$ given in Definition~\ref{def:Tsigmatheta}, and we take $p_1$, $p_2\in\Q$ such that $J=[p_1,p_2]\subset I$, then $y(\cdot)$ does not necessarily belong to $\K_j^J$. Thus, we need to consider the next technical lemma that will be used later.
\begin{lem}\label{lem:conv-subint}
Let $\Theta$ be a suitable set of moduli of continuity as in {\rm Definition~\ref{def:modcont}}.
Let $\{\nu_n\}_{n\geq 1}$ be a sequence in $\PM$ that converges to $\nu_0\in \PM$ in the topology $\sigma_\Theta$. Then, for each integer $j\geq 1$ and each interval $I=[q_1,q_2]$, $q_1,q_2\in\Q$,
\[
\lim_{n\to\infty}\sup_{y(\cdot)\in\K_j^I}\left| \int_{p_1}^{p_2}\left[d(\nu_n)_{y(s)}- d(\nu_0)_{y(s)}\right]\right|=0
\]
whenever $p_1$, $p_2\in \Q$ and $q_1\le p_1 < p_2\le q_2$.
\end{lem}
\begin{proof}
For each $y(\cdot)\in\K_j^I$, we define $\wt y\colon I\to B_j$ by
$\wt y(t)=y(p_1)$ if $t\in [q_1,p_1]$, $\wt y(t)=y(t)$ if $t\in [p_1,p_2]$, and $\wt y(t)=y(p_2)$ if $t\in [p_2,q_2]$. The map $\wt y(\cdot)\in\K_j^I$ and
\[
 \left| \int_{p_1}^{p_2}\left[d(\nu_n)_{y(s)}- d(\nu_0)_{y(s)}\right]\right| =  \left| \int_{p_1}^{p_2}\left[d(\nu_n)_{\wt y(s)}- d(\nu_0)_{\wt y(s)}\right]\right| \leq I_n^1+ I_n^2+I_n^3 \,,
\]
where  $ I_n^1=\big| \int_{q_1}^{q_2}\left[d(\nu_n)_{\wt y(s)}- d(\nu_0)_{\wt y(s)}\right]\big|$, $ I_n^2=\big| \int_{q_1}^{p_1}\left[d(\nu_n)_{y(p_1)}- d(\nu_0)_{ y(p_1)}\right]\big|$ and
 $ I_n^3=\big| \int_{p_2}^{q_2}\left[d(\nu_n)_{y(p_2)}- d(\nu_0)_{y(p_2)}\right]\big|$, from which the result can easily be deduced because $\nu_n\xrightarrow[]{\sigma_\Theta}\nu_0$, $\wt y(\cdot)\in \K_j^I$ and the constant functions $y(p_1)$ and $y(p_2)$ belong  to $\K_j^{[q_1,p_1]}$ and  $\K_j^{[p_2,q_2]}$,  respectively.
 \end{proof}
Next, we introduce some topological properties of the $m$-bounds and $l$-bounds depending on the standard definitions for families of positive measures in $\M_c^+$.
\begin{defn}\label{def:families}
A family of positive measures $\{\mu_\alpha\mid \alpha\in\Gamma\}\subset \M_c^+$ 
\begin{itemize}
\item[(i)] is \emph{bounded} if there is a positive constant $c>0$ such that $\mu_\alpha[t,t+1]\leq c$ for each $t\in\R$  and~$\alpha\in\Gamma$;
\item[(ii)] is \emph{equicontinuous}
if for each $\ep>0$ there is a $\delta(\ep)>0$ such that $\mu_\alpha[s,t]<\ep$ whenever $0\leq t-s < \delta(\ep)$ and $\alpha\in\Gamma$;
\item[(iii)] is \emph{uniformly integrable} if given $\ep>0$ there is a $\delta(\ep)>0$ such that if $A$ is a Borel set satisfying $A\subset[t,t+1]$ for some $t\in\R$ and $l(A)< \delta(\ep)$, then $\mu_\alpha(A)< \ep$ for each $\alpha\in \Gamma$.
\end{itemize}
\end{defn}
\begin{rmk}\label{rmk:implications}
According to the previous definitions, uniform integrability implies equicontinuity, and equicontinuity implies boundedness. Moreover, notice that condition (ii) allows for singular continuous measures with respect to the Lebesgue measure, whereas with condition (iii) all the measures of the family are absolutely continuous with respect to the Lebesgue measure.
\end{rmk}

\begin{defn}
A family   $\{\nu_\alpha \mid \alpha\in\Gamma\}\subset\PM$ of parametric b-measures with $m$\nbd-bounds $\{\{m_j^\alpha\}_{j\geq 1}\mid  \alpha\in\Gamma\}\subset \M_c^+$ and $l$-bounds
$\{\{l_j^\alpha\}_{j\geq 1}\mid  \alpha\in\Gamma\}\subset \M_c^+$
\begin{itemize}
\item[(i)] \emph{has bounded $m$-bounds} (resp., \emph{bounded $l$-bounds}) if the set $\{m_j^\alpha\mid \alpha\in\Gamma\}$ (resp.,  $\{l_j^\alpha\mid \alpha\in\Gamma\}$) is bounded  for each $j\geq 1$;
\item[(ii)] \emph{has equicontinuous $m$-bounds} (resp., \emph{equicontinuous $l$-bounds}) if the set $\{m_j^\alpha\mid \alpha\in\Gamma\}$ (resp., $\{l_j^\alpha\mid \alpha\in\Gamma\}$) is equicontinuous for each $j\geq 1$;
\item[(iii)] \emph{has uniformly integrable $m$-bounds} (resp., \emph{uniformly integrable $l$-bounds}) if  $\{m_j^\alpha\mid \alpha\in\Gamma\}$ (resp., $\{l_j^\alpha\mid \alpha\in\Gamma\}$) is uniformly integrable for each $j\geq 1$.
\end{itemize}
\end{defn}
Since in this work we are not going to require the property of uniform integrability, we will only prove that the properties of boundedness and equicontinuity of the $m$-bounds/$l$-bounds of a set  $E\subset \PM$  are inherited by its closure $\overline E=\cls_{(\PM,\T)}(E)$ for $\T\in\{\sigma_\Theta, \sigma_D\}$, but the same happens with the property of uniform integrability. Specifically, $\mathrm{cls}_{(\PM,\T)}(E)$ represents the closure in $(\PM,\T)$ of the subset of parametric b-measures  $E$ with the induced topology.  We recall the following result for the vague topology $\widetilde{\sigma}$ in $\M^+$, which will be used in the proof.

\begin{prop}\label{prop:boundedness}
$($\cite[Theorem~15.7.5, p.170]{book:Kall}$)$ A subset $M_0\subset\M^+$ is relatively compact in the vague topology $\wt\sigma$ if and only if $\sup_{\mu\in M_0}\mu (B)<\infty$ for every bounded Borel set $B\in\B$.
\end{prop}

\begin{prop}\label{prop:EclsE}
Let $\T$ be one of the topologies $\{\sigma_\Theta, \sigma_D\}$ defined above.
\begin{itemize}
\item[(i)] If a subset $E\subset \PM$ has bounded $m$-bounds (resp., bounded $l$-bounds), then $\overline E=\cls_{(\PM,\T)}(E)$ also has bounded $m$-bounds (resp., bounded $l$-bounds).
\item[(ii)] If a subset $E\subset \PM$ has equicontinuous $m$-bounds (resp., equicontinuous $l$-bounds), then $\overline E=\cls_{(\PM,\T)}(E)$ also has equicontinuous $m$-bounds (resp., equicontinuous $l$-bounds).
 \end{itemize}
\end{prop}
\begin{proof}
Since $\sigma_D\leq \sigma_\Theta$, it is enough to prove the result for $\sigma_D$ .

(i)
Let $E=\{\nu_\alpha \mid \alpha\in\Gamma\}\subset\PM$ be a subset of parametric b-measures with bounded $m$\nbd-bounds $\{m_j^\alpha\mid \alpha\in\Gamma\}\subset \M_c^+$ for each $j\geq 1$. Given $\mu\in\cls_{(\PM,\sigma_D)}(E)$, there is a sequence $\{\nu_{\alpha_k}\}_{k\geq 1}\subset E$  such that $\mu=\lim_{k\to\infty}\nu_{\alpha_k}$  in $\sigma_D$, i.e., for every $I=[q_1,q_2]$, $q_1,q_2\in\Q$ and every $y\in D$, $\lim_{k\to\infty} |(\nu_{\alpha_k})_y(I)-\mu_y(I)|=0$. Since in particular $\mu\in \PM$, there exist  $\{m_j^\mu, l_j^\mu\}_{j\geq 1}\subset \M_c^+$ $m$-bounds, $l$-bounds, respectively, of $\mu$ on $B_j$. The boundedness of $m_j^\mu$ cannot be studied directly. Thus, we build new $m$-bounds for $\mu$ by using the $m$-bounds of the sequence $\{\nu_{\alpha_k}\}_{k\geq 1}\subset E$.

Fixed $j\geq 1$, from Proposition~\ref{prop:boundedness}, taking a subsequence if necessary,  $\{m_j^{\alpha_{k}}\}_{ k\geq 1}$  converges to some $m_j$ in $(\M^+,\widetilde{\sigma})$. We claim that the continuous part of $m_j$ is an $m$-bound for $\mu$ on $B_j$.   First, if $I=[q_1,q_2]$, $q_1,q_2\in\Q$ and $y\in B_j\cap D$,
\begin{equation}\label{eq:aux 2}
|\mu_y(I)|=\lim_{k\to\infty} |(\nu_{\alpha_k})_y(I)| \quad\mbox{and}\quad |(\nu_{\alpha_k})_y(I)|\leq |(\nu_{\alpha_k})_y|(I) \leq m_j^{\alpha_k}(I)\,.
\end{equation}
By the regularity of the measures  $m_j^{\alpha_k}$, and the vague convergence, respectively,
\begin{align*}
 & m_j^{\alpha_k}(I)  =\inf \bigg\{ \int_\R \phi(s)\,dm_j^{\alpha_k}\, \Big|\,\,  \phi\in C_c^+(\R)\,, \;\phi|_{I}\equiv 1  \bigg\}\,, \\
 & \lim_{k\to\infty} \int_\R \phi(s)\,dm_j^{\alpha_{k}}= \int_\R \phi(s)\,dm_j \;\text{ for each } \phi\in C_c(\R)\,.
  \end{align*}
Since also $m_j$ is a regular measure, we can deduce from \eqref{eq:aux 2} and the previous formulas that $|\mu_y(I)|\leq m_j(I)$. We note that we can remove the purely discontinuous part, if any, from the measure $m_j$. That is, if $m_j=m_{j,c} + m_{j,pd}\in \M_c^+\oplus\M^+_{pd}$, then $|\mu_y(I)|\leq m_{j,c}(I)$ for every $I=[q_1,q_2]$, $q_1,q_2\in\Q$ and $y\in B_j\cap D$, because $\mu_y\in \M_c(\B)$ and $m_{j,pd}$ is supported on a countable set $\mathcal{N}$ with $\mu_y(I\cap\mathcal{N})=0$.

Now, approximating a compact interval $J$  by intervals with rational endpoints and using the regularity of $\mu_y$  we deduce that
$| \mu_y (J)|\leq m_{j,c} (J)$.  Notice that since $\mu_y\in \M_c(\B)$ and $m_{j,c}\in \M_c^+$  do not have purely discontinuous part, the inequality follows also for each open interval $J$. Since any bounded open set $U$  is a disjoint countable union of bounded open subintervals, the inequality holds for $U$.
Again, by the regularity property, the inequality  $|\mu_y (A)|\leq m_{j,c} (A)$ follows for every bounded Borel set $A$ and each  $y\in B_j\cap D$. And, if $y\in B_j\setminus D$, then $|\mu_y (A)- \mu_x(A)|\leq |\mu_y - \mu_x|(A) \leq l_j^\mu(A) \,\omega_j(|y-x|)$ for all $A\in\B$, so we only need to approximate $y$ by points $x\in B_j\cap D$ to extend the result to all $y\in B_j$. Finally, by the definition of $|\mu_y|$ in \eqref{eq:total var} and since $m_{j,c}$ is a measure, it follows that $|\mu_y|(A)\leq m_{j,c}(A)$ for all $A\in\B$ and $y\in B_j$, as claimed.

To finish,  if $m_j^\alpha[t,t+1]\leq c_j$ for $t\in\R$  and~$\alpha\in\Gamma$, for the $m$-bound $m_{j,c}$ on $B_j$, taking for each $t\in\R$ a map $\phi\in C_c(\R)$, which depends on $t$ and satisfies $0\leq \phi(s)\leq 1$ for each $s\in\R$, $\phi(s)=1$ for $s\in [t,t+1]$ and $\phi(s)=0$ for $s\in(-\infty,t-1]\cup [t+2,\infty)$, and noting that $m_{j,c}\leq m_j$, we deduce that for every $t\in\R$,
\begin{equation}\label{eq:aux 3}
  m_{j,c}[t, t+1]\leq \int_\R\phi(s) \,dm_{j}=\lim_{k\to\infty}\int_\R\phi(s)\, dm_j^{\alpha_{k}}\leq \sup_{\alpha\in\Gamma} m_j^\alpha[t-1,t+2]<3\,c_j\,.
\end{equation}
Since the same arguments apply to each $\mu\in \overline E$, and each $j\geq 1$,  we can conclude that  $\overline E=\cls_{(\PM,\sigma_D)}(E)$ has bounded $m$-bounds, as claimed.  We omit the proof in the case of bounded $l$-bounds, because it is analogous.
\par\smallskip
(ii) We maintain the notation used in (i) and assume that for each $j\geq 1$, the set of $m$-bounds $\{m_j^\alpha\mid \alpha\in\Gamma\}\subset \M_c^+$ is equicontinuous, i.e., given $\ep>0$ there is a $\delta(\ep)>0$ such that $m_j^\alpha[s,t]<\ep$ whenever $0\leq t-s<\delta(\ep)$, for all $\alpha\in\Gamma$. Since equicontinuity implies boundedness, for each $\mu \in\overline{E}$ we can consider the $m$-bound $m_{j,c}$ built in (i). Then, similar arguments as those in \eqref{eq:aux 3} applied to $m_{j,c}[s, t]$ allow us to conclude that $m_{j,c}[s, t]<3\,\ep$ whenever $0\leq t-s<\delta(\ep)$. Thus,  $\overline{E}$ has equicontinuous $m$-bounds. The case of equicontinuous $l$-bounds
is also omitted by similarity of reasoning. The proof is finished.
\end{proof}
\par
Next, we want to show the equivalence of the introduced topologies under the existence of bounded $l$-bounds,  an assumption that is satisfied in many applications.
Before that, we give a technical result.
\begin{lem}\label{lem:uniformtau}
Let $E\subset \PM$ have bounded $l$-bounds $\{\{l_j^\nu\}_{j\geq 1}\mid  \nu\in E\}$ and let $\Theta$ be a suitable set of moduli of continuity as in {\rm Definition~\ref{def:modcont}}. Fixed  $I=[q_1,q_2]$ $(q_1,q_2\in\Q)$ and an integer $j\geq 1$, consider the set $\K_j^I$ in \eqref{eq:KjI}. Then, given $\ep>0$,  there is a $\delta(\ep)>0$ such that for each tagged $\delta(\ep)$\nbd-partition $\Delta^{\delta(\ep)}=\{(\tau_i, [t_i,t_{i+1}])\mid i=1,\ldots, k-1\}$ of $I$, for all $y(\cdot)\in \K_j^I$ and $\nu\in E$ it holds:
\begin{equation*}
		\left|\int_I d\nu_{y(s)}-S_\nu^{y(\cdot)}\big(\Delta^{\delta(\ep)}\big)\Bigg|=\Bigg| \int_I d\nu_{y(s)}-\sum_{i=1}^{k-1} \int_{t_i}^{t_{i+1}}d\nu_{y(\tau_i)}\right|<\ep\,.
	\end{equation*}
\end{lem}	
\begin{proof}
First of all, we affirm that the statement of Lemma~\ref{lem:tau} holds uniformly for all  $y(\cdot)\in \K_j^I$ and $\nu\in E$. Note that the functions $y(\cdot)\in \K_j^I$ are uniformly bounded by $j$ and uniformly equicontinuous and, moreover, we can take an $M>0$ such that $l_j^\nu (I)\leq M$ for all $\nu\in E$ because $E$ has bounded $l$-bounds.

Then, arguing as in the proof of Proposition~\ref{prop:defIntegral} it is easy to prove that given $\ep>0$,  there is a $\delta>0$ such that $\big|S_\nu^{y(\cdot)}(\Delta_1^\delta)- S_\nu^{y(\cdot)}(\Delta_2^\delta)\big|<\ep$
for all tagged $\delta$-partitions $\Delta_1^\delta$  and $\Delta_2^\delta$ of $I$,  uniformly for $y(\cdot)\in \K_j^I $ and  $\nu\in E$. That is, the Riemann sums satisfy the Cauchy property uniformly for all $y(\cdot)\in \K_j^I$ and $\nu\in E$, and this is equivalent to what we wanted to prove. The proof is finished.
\end{proof}
\begin{thm}\label{thm:equivalenceWeakTopologies}
	Let $E$ be a subset of $\PM$ with bounded $l$-bounds. Consider a dense and countable subset $D$ of $\R^N$ and a suitable set of moduli of continuity $\Theta$ as in {\rm Definition~\ref{def:modcont}}. Then,
\begin{equation*}	(E,\sigma_D)=(E,\sigma_\Theta)\qquad\text{and}\qquad\mathrm{cls}_{(\PM,\sigma_D)}(E)=\mathrm{cls}_{(\PM,\sigma_\Theta)}(E)\, .
	\end{equation*}
\end{thm}		
\begin{proof}
As noted in Remark \ref{rmk:orden top}, $\sigma_D\leq \sigma_\Theta$. Thus, it suffices to prove that, if a sequence $\{\nu_n\}_{n\geq 1}\subset E$ converges to $\nu_0\in \PM$  in the topology $\sigma_D$, that is,
\begin{equation}\label{convD}
	\lim_{n\to\infty}\left|\,\int_I d(\nu_n)_y-\int_I d(\nu_0)_y\,\right|=0\quad \text{for all } y\in D\,,\, I=[q_1,q_2]\,\, (q_1,\, q_2\in \Q)\,,
\end{equation}
then, it also converges to $\nu_0$ in $\sigma_\Theta$, i.e., for all $j\geq 1$ and $I=[q_1,q_2]\,(q_1,q_2\in \Q)$,
\begin{equation*}\label{convTheta}
\lim_{n\to\infty}\sup_{y(\cdot)\in \K_j^I }\left|\,\int_I d(\nu_n)_{y(s)}- \int_I d(\nu_0)_{y(s)}\,\right|=0\,.
\end{equation*}
The main ideas come from the proof of~\cite[Theorem~2.33]{tesis:Iacopo} for Carath\'{e}odory functions, but the adaptation to the case of parametric b-measures requires an explanation of some nontrivial technical details.  Fix a compact interval $I=[q_1,q_2]$ with $q_1$, $q_2\in \Q$ and $j\geq 1$, and denote by $l_j^n \in\mathcal M_{c}^+$ an $l$-bound of $\nu_n$ on $B_j$, for each $n\geq 0$. By Proposition~\ref{prop:EclsE}, the set $\mathrm{cls}_{(\PM,\sigma_D)}(E)$ has bounded $l$-bounds. Then, there is a constant $c>0$ such that $l_j^n(I)\leq c$ for each $n\geq 0$.

Now, fix $\ep>0$ and take $\rho>0$ so that $\omega_j(\rho)<\ep/(6\,c)$. Since $B_j\subset \R^N$ is compact, and $D$ is dense in $\R^N$, there are $y_1,y_2,\ldots,y_{\ell_0}\in B_j\cap D$ such that $B_j\subset\bigcup_{\ell=1}^{\ell_0} \mathring{B}_\rho(y_\ell)$, where $\mathring{B}_\rho(y)$ denotes the open ball of $\R^N$ of radius $\rho$ centered at $y\in\R^N$. Subordinate  to the former open cover, there is a continuous partition of  unity, that is, there exist continuous functions $\phi_\ell\colon \R^N \mapsto [0,1]$ for $\ell=1,\ldots,\ell_0$, which in fact can be taken to be of class $C^1$, such that
\begin{equation}\label{particion}
\supp(\phi_\ell)\subset \mathring{B}_\rho(y_\ell)\, \qquad \mathrm{and}\qquad \sum_{\ell=1}^{\ell_0} \phi_\ell(y)=1 \quad \forall\, y\in B_j \, .
\end{equation}
From Proposition~\ref{prop:phinu}, the maps $\nu_n^*\colon\R^N\to \M_c(\B)$, $y\mapsto (\nu_n^*)_y$ given by
\[
(\nu_n^*)_y=\sum_{\ell=1}^{\ell_0} \phi_\ell(y)(\nu_n)_{y_\ell} \quad \text{ for each }y\in \R^N \text{ and } n\geq 0
\]
permit us to define a family of parametric b-measures $\{\nu_n^*\mid n\geq 0\}$ with respect to the trivial family of moduli of continuity given by the identity map.

Now, for each $y(\cdot)\in \K_j^I$, and each $n\geq 1$, we consider a tagged $\delta$-partition of $I$,  $\Delta^\delta=\{(\tau_i,[t_i,t_{i+1}])\mid i=1,\ldots, k-1\}$ with rational points $t_i$, for convenience, and apply the triangle inequality to write
\begin{align}\label{3integrals}
\left|\,\int_I d(\nu_n)_{y(s)}-\int_I    d(\nu_0)_{y(s)}\,\right| \leq &  \left|\,\int_I d(\nu_n)_{y(s)} - S_{\nu_n}^{y(\cdot)}(\Delta^\delta)\right| +
\left| S_{\nu_n}^{y(\cdot)}(\Delta^\delta) - S_{\nu_n^*}^{y(\cdot)}(\Delta^\delta)\right|   \nonumber \\ &
  + \left| S_{\nu_n^*}^{y(\cdot)}(\Delta^\delta) - S_{\nu_0^*}^{y(\cdot)}(\Delta^\delta)\right|   + \left|S_{\nu_0^*}^{y(\cdot)}(\Delta^\delta) -S_{\nu_0}^{y(\cdot)}(\Delta^\delta)\right|   \nonumber \\ & + \left| S_{\nu_0}^{y(\cdot)}(\Delta^\delta) - \int_I d(\nu_0)_{y(s)} \right| .
 \end{align}
From Lemma~\ref{lem:uniformtau} applied to $\mathrm{cls}_{(\PM,\sigma_D)}(E)$, associated with the fixed $\ep>0$, there exists a $\delta>0$ such that for each tagged $\delta$\nbd-partition $\Delta^\delta$ of $I$,
\begin{equation*}
a_n:=\left|\,\int_I d(\nu_n)_{y(s)}- S_{\nu_n}^{y(\cdot)}(\Delta^\delta)\right|<\frac{\ep}{6}\,\quad\hbox{for all } n\geq 0 \;\,\hbox{and } y(\cdot)\in \K_j^I\,.
\end{equation*}
This applies to the first and last terms in the sum in \eqref{3integrals}. For the second and fourth terms,  since for every $y(\cdot)\in \K_j^I$, $y(\tau_i)\in B_j$ for all $i=1,\ldots, k-1$, applying~\eqref{particion} and Proposition~\ref{prop:Phinuintegral}, and denoting $I_i=[t_i,t_{i+1}]$, we write for each $n\geq 0$,
\begin{align*}
b_n:=	\left|  S_{\nu_n}^{y(\cdot)}(\Delta^\delta)- S_{\nu^*_n}^{y(\cdot)}(\Delta^\delta)\right|& = \left| \, \sum_{i=1}^{k-1}\int_{I_i} d(\nu_n)_{y(\tau_i)}-\sum_{i=1}^{k-1} \int_{I_i} \sum_{\ell=1}^{\ell_0}\phi_\ell(y(\tau_i))\,d(\nu_n)_{y_\ell}\right| \\
	& =\left|\,  \sum_{i=1}^{k-1}\int_{I_i} \sum_{\ell=1}^{\ell_0}\phi_\ell(y(\tau_i)) \left(d(\nu_n)_{y(\tau_i)}-d(\nu_n)_{y_\ell}\right)\right|.
\end{align*}
 Therefore, from the bounded $l$-bounds of $\mathrm{cls}_{(\PM,\sigma_D)}(E)$ and again using~\eqref{particion}, we obtain that for each $n\geq 0$ and $y(\cdot)\in \K_j^I$,
\begin{align*}
b_n& \leq
	\sum_{i=1}^{k-1}\left[ \sum_{\ell=1}^{\ell_0}\phi_\ell(y(\tau_i)) \left|(\nu_n)_{y(\tau_i)}-(\nu_n)_{y_\ell}\right| (I_i)\right]\\
	&\leq \sum_{i=1}^{k-1} \left[\sum_{\ell=1}^{\ell_0}\phi_\ell(y(\tau_i))\,\omega_j(|y(\tau_i)-y_\ell|)\right]l_j^n(I_i)
 \leq \omega_j(\rho) \sum_{i=1}^{k-1} l_j^n(I_i)< \frac{\ep}{6\,c}\, l_j^n(I)\leq\frac{\ep}{6}\,.
	\end{align*}
Finally, for the third term in the sum in \eqref{3integrals} we have
\[
c_n:= \left| S_{\nu_n^*}^{y(\cdot)}(\Delta^\delta) - S_{\nu_0^*}^{y(\cdot)}(\Delta^\delta)\right| \leq \sum_{i=1}^{k-1} \left[ \sum_{\ell=1}^{\ell_0}\phi_\ell(y(\tau_i))\left|\int_{I_i} \left[d(\nu_n)_{y_\ell}-d(\nu_{0})_{y_\ell}\right]\right| \right]\,.
\]
Since $\{\nu_n\}_{n\geq 1}\subset E$ converges  to $\nu_0$ in $\sigma_D$, that is, \eqref{convD} holds, there is an $n_0$ such that for each $n\geq n_0$, each $\ell=1,\ldots,\ell_0$ and  each $i=1,\ldots, k-1$,
\[
\left| \int_{I_i}\left[d(\nu_n)_{y_\ell}-d(\nu_{0})_{y_\ell} \right]  \right|< \frac{\ep}{3\, (k-1)}\,.
\]
Then, using~\eqref{particion} once more, $c_n<\ep/3$ for all $y(\cdot)\in\K^I_j$ and $n\geq n_0$. Collecting the previous facts,  we can conclude that $\{\nu_n\}_{n\geq 1}$ converges  to $\nu_0$ in $\sigma_\Theta$, as we wanted. The proof is finished.
\end{proof}
We conclude this section with the crucial result of compactness for sets of parametric b-measures with equicontinuous $m$-bounds and bounded $l$-bounds.
\begin{thm}\label{th-compacidad}
Let $E\subset\PM$ have equicontinuous $m$-bounds and bounded $l$-bounds and let $\T\in\{\sigma_D, \sigma_\Theta\}$. Then, $\mathrm{cls}_{(\PM,\T)}(E)$ is compact.
\end{thm}
\begin{proof}
By Theorem~\ref{thm:equivalenceWeakTopologies}, with the hypothesis of bounded $l$-bounds the topologies $\sigma_D$ and $\sigma_\Theta$ are equivalent on $E$,  
so that it suffices to prove the result for the topology $\sigma_D$, for a countable and dense subset $D\subset\R^N$.
Let us prove that, given a sequence $\{\nu_n\}_{n\geq 1}\subset E$ there exists a subsequence $\{\nu_{n_k}\}_{k\geq 1}$ which converges in the topology $\sigma_D$ to some $\nu\in\PM$.
For each $n\geq 1$ and $j\geq 1$, respectively denote by $m_j^n$ and $l_j^n$ an $m$-bound and an $l$-bound of $\nu_n$ on $B_j$. Now, for each $n\geq 1$, define the map
\[
F_n\colon \R\times\R^N \to\R\,,\quad (t,y)\mapsto F_n(t,y):= \int_0^t d(\nu_n)_y\,.
\]
With this definition, for $n\geq 1$ and for each $y\in\R^N$, $(\nu_n)_y[t_1,t_2]=F_n(t_2,y)-F_n(t_1,y)$ for each real interval $[t_1,t_2]$ and
$F_n(0,y)=0$ because $(\nu_n)_y\in \M_c(\B)$.

Let us see that the family $\mathcal{F}=\{F_n\mid n\geq 1\}$ is uniformly equicontinuous on the compact set $[-j,j]\times B_j$, for each $j\geq 1$. Take $t_1,t_2\in [-j,j]$ with $t_1<t_2$ and $y_1,y_2\in B_j$ and write for $n\geq 1$, assuming that $t_1>0$ without loss of generality,
\begin{align*}
    |&F_n(t_2,y_2)  - F_n(t_1,y_1)| = \left| \int_0^{t_2} d(\nu_n)_{y_2}- \int_0^{t_1} d(\nu_n)_{y_1} \right| \\
    &\; \leq  \left| \int_0^{t_2} d(\nu_n)_{y_2}- \int_0^{t_1} d(\nu_n)_{y_2} \right| +  \left| \int_0^{t_1} d(\nu_n)_{y_2}- \int_0^{t_1} d(\nu_n)_{y_1} \right|   \\[1pt]
    & \;\leq  |(\nu_n)_{y_2}| [t_1,t_2] + |(\nu_n)_{y_2}-(\nu_n)_{y_1}|[0,t_1]
    \leq   m_j^n [t_1,t_2] + l_j^n[-j,j]\,\omega_j(|y_2-y_1|)\,.
\end{align*}
From the hypotheses of equicontinuous $m$-bounds and bounded $l$-bounds for the set $E$, the uniform equicontinuity follows directly. Furthermore, the family
$\mathcal{F}$ is bounded at every fixed $(t,y)\in \R\times\R^N$: take $j\geq 1$ such that $y\in B_j$, and assuming that $t>0$ without loss of generality, $|F_n(t,y)|\leq  |(\nu_n)_{y}|[0,t]\leq  m_j^n[0,t]$ for all $n\geq 1$, and the $m$-bounds are, in particular, bounded (see Remark~\ref{rmk:implications}).

Since $\R\times\R^N$ is separable, by the Arzel\`{a}-Ascoli theorem and a diagonalization process there exists a subsequence $\{F_{n_k}\}_{k\geq 1}$ that is uniformly convergent on every compact subset of $\R\times\R^N$. Hence, there is a continuous function $F\colon \R\times\R^N \to\R$ such that $F_{n_k}\to F$ as $k\to\infty$ uniformly on compact sets.

In order to build a parametric b-measure from the continuous map $F(t,y)$, let us prove that for each  $y\in \R^N$ fixed, the map $F(\cdot,y)$ is of bounded variation on every compact set $[a,b]$. First, because $E$ has bounded $m$-bounds, it is easy to check that the family $\{F_n(\cdot,y)\mid n\geq 1\}$ is uniformly of bounded variation on $[a,b]$. To see it, take $j\geq 1$ such that $y\in B_j$ and $c_0>0$ such that
$m_j^n [a,b]\leq c_0$ for every $n\geq 1$. Then, given any integer $\ell\geq 1$ and any partition $a=t_0<t_1<\ldots<t_\ell=b$,
\begin{equation}\label{eq:boundedvar}
\begin{split}
\sum_{i=1}^\ell \,|F_n(t_i,y)-F_n(t_{i-1},y)| & = \sum_{i=1}^\ell |(\nu_n)_y[t_{i-1},t_i]|\leq \sum_{i=1}^\ell|(\nu_n)_y|[t_{i-1},t_i] \\
 &\leq \sum_{i=1}^\ell m_j^n[t_{i-1},t_i] = m_j^n [a,b]\leq c_0\,
\end{split}
\end{equation}
for all $n\geq 1$. Second, for the map $F$, for each $\ell\geq 1$ we can take $k=k(\ell)$ big enough so that $|F(t,y)-F_{n_k}(t,y)|<1/\ell$ for all $t\in [a,b]$, and then  $\sum_{i=1}^\ell |F(t_i,y)-F_{n_k}(t_{i},y)|\leq 1$ and $\sum_{i=1}^\ell |F_{n_k}(t_{i-1},y)-F(t_{i-1},y)|\leq 1$. Thus, for any $\ell\geq 1$ and any partition of $[a,b]$, applying the triangle inequality and \eqref{eq:boundedvar} for the map $F_{n_k}$,  we get
\[
\sum_{i=1}^\ell |F(t_i,y)-F(t_{i-1},y)| \leq 2 +c_0\,.
\]
As a consequence, $F(\cdot,y)$ is of bounded variation on every compact interval $[a,b]$.

Therefore, we can consider $\nu_y$ the Lebesgue-Stieltjes real measure on $[a,b]$  associated with $F(\cdot,y)$, which is uniquely determined and satisfies that $\nu_y[t_1,t_2]=F(t_2,y)-F(t_1,y)$ for intervals in  $[a,b]$, by construction. Thus, by considering growing intervals $[-j,j]$ one can extend $\nu_y$ to the ring of bounded Borel sets in $\R$. Note that, since $F(\cdot,y)$ is continuous, $\nu_y(\{t\})=0$ for every $t\in\R$ and we can affirm that $\nu_y\in \M_c(\B)$. Furthermore, $\nu\colon\R^N\to \M_c(\B)$ defined in this way is the limit of $\nu_{n_k}$ in the topology $\sigma_D$. In fact,  for every $y\in \R^N$ and every interval $I=[q_1,q_2]$,
\begin{equation}\label{eq:aux}
\left|\int_I d(\nu_{n_k})_y - \int_I d\nu_y\right| = \big|F_{n_k}(q_2,y)-F_{n_k}(q_1,y) - (F(q_2,y)-F(q_1,y))\big| \to 0
\end{equation}
as $k\to\infty$. That is, $\nu_y(I)= \lim_{k\to\infty} (\nu_{n_k})_y(I)$ for all $y\in \R^N$ and every interval $I\subset\R$. To complete the proof, it remains to build $m$-bounds and $l$-bounds for~$\nu$. As in the proof of Proposition~\ref{prop:EclsE}, we can do this by using the $m$-bounds and $l$\nbd-bounds of $\nu_{n_k}$. Namely, fixed $j\geq 1$, from Proposition~\ref{prop:boundedness}, taking a subsequence if necessary,
$\{m_j^{n_{k}}\}_{ k\geq 1}$  converges to some $m_j$  and $\{l_j^{n_{k}}\}_{ k\geq 1}$ converges to some $l_j$ in $(\M^+,\widetilde{\sigma})$. Since \eqref{eq:aux} holds for all $y\in \R^N$ and every interval $I\subset\R$, it is now easier than in the mentioned proof to check that the continuous part of $m_j$ (resp., of $l_j$) is an $m$-bound (resp., an $l$-bound) of $\nu$ on $B_j$.   Thus, $\nu\in\PM$ and the proof is finished.
\end{proof}
\section{Continuity of time translations and the hull of a parametric  b-measure}\label{sec hull}
As in the previous section, we consider the set  $\PM$ of parametric b-measures on $\R^N$ with respect to a fixed family of moduli of continuity  $\{\omega_j\}_{j\geq 1}$ (see Definition~\ref{defi:parameasure}). In order to build the hull of a parametric b-measure, we define the translation of a parametric b-measure by a number $t\in\R$.
Note that, given a measure $\mu\in\M^+$ (resp., a b-measure $\mu\in \M_c(\B)$) and a real number $t$, $\mu{\cdot}t$ defined by $\mu{\cdot}t(A) := \mu(A+t)$ for each $A\in \mathcal{A}$ (resp., $A\in \B$) is a measure (resp., a b-measure).
\begin{lem}\label{lem:translmeaseure}
Given a parametric b-measure $\nu\in\PM$ with $m$\nbd-bounds $\{m_j\}_{j\geq 1}\subset \M_c^+$ and $l$-bounds
$\{l_j\}_{j\geq 1}\subset \M_c^+$,  and a real number $t\in\R$, we consider the map $\nu{\cdot}t\colon \R^N\to \M_c(\B),\;y\mapsto (\nu{\cdot}t)_y:= \nu_y{\cdot}t$. Then, $\nu{\cdot}t$ is a parametric b-measure  with $m$-bounds $\{m_j{\cdot} t\}_{j\geq 1}$ and $l$-bounds $\{l_j{\cdot}t\}_{j\geq 1}$.
\end{lem}
\begin{proof}
Fixed an integer $j\geq 1$, from the definition of the total variation of a b-measure in \eqref{eq:total var}, we deduce that   for all $A\in \B$ and $y, y_1, y_2\in B_j\subset\R^N$,
\begin{align*}
   |(\nu{\cdot}t)_y|(A)& =|\nu_y|(A+t)\leq m_j(A+t)=m_j{\cdot}t(A)\,,\\
   |(\nu{\cdot}t)_{y_1}-(\nu{\cdot}t)_{y_2}|(A)& =|\nu_{y_1}-\nu_{y_2}|(A+t) \leq l_j(A+t)\,\omega_j(|y_1-y_2|)\\& =l_j{\cdot}t(A)\,\omega_j(|y_1-y_2|)\,.
\end{align*}
With this, the proof is finished.
\end{proof}
We give a formula for the integration of a translated parametric b-measure along a continuous map $y\colon [a,b]\to\R^N$ to be used in the proof of Theorem~\ref{conti-time-trns}.
\begin{lem}\label{lem:integralnut}
Let $\nu\in\PM$ be a parametric b-measure on $\R^N$, $t\in\R$ a real number, and  $y\colon [a,b]\to \R^N$ a continuous map. Then,
\[
\int_a^b d(\nu{\cdot}t)_{y(s)}=\int_{a+t}^{b+t} d\nu_{y(s-t)}\,.
\]
\end{lem}
\begin{proof} It is enough to notice that given
$ \Delta^\delta =\{(\tau_i, [t_i,t_{i+1}])\mid i=1,\ldots, k-1\}$ any $\delta$-partition of $[a,b]$, $ \Delta^\delta_t =\{(\tau_i+t, [t_i+t,t_{i+1}+t])\mid i=1,\ldots, k-1\}$ is a $\delta$-partition of $[a+t,b+t]$ and the corresponding sums
\[ S_{\nu{\cdot}t}^{y(\cdot)}\big(\Delta^\delta \big)=\sum_{i=1}^{k-1} \int_{t_i}^{t_{i+1}}d(\nu{\cdot}t)_{y(\tau_i)} \quad \text{and} \quad S_{\nu}^{y(\,\cdot\, -t)}\big(\Delta_t^\delta \big)=\sum_{i=1}^{k-1} \int_{t_i+t}^{t_{i+1}+t}d\nu_{y(\tau_i+t -t)}\]
coincide.
\end{proof}
Next, we prove the continuity of time translations on sets with equicontinuous $m$-bounds. This result is the counterpart of Theorem~3.1 in~\cite{paper:LNO2} for the translation of Carath\'{e}odory functions.
\begin{thm}\label{conti-time-trns}
Let $\Theta$ be a suitable set of moduli of continuity as in {\rm Definition~\ref{def:modcont}},  $E$ be a subset of $\PM$ with equicontinuous $m$-bounds and $\overline E=\cls_{(\PM,\sigma_\Theta)}(E)$. Then, the map
\begin{equation*}
\R\times \overline E\to \PM\, ,\qquad (t,\nu)\mapsto \nu{\cdot}t\,
\end{equation*}
is continuous.
\end{thm}
\begin{proof}
Notice that the map is well defined thanks to Lemma~\ref{lem:translmeaseure}. Let $\{\nu_n\}_{n\geq 1}\subset \overline E$ be a sequence that converges to $\nu_0$ in $(\PM,\sigma_\Theta)$ and $\{t_n\}_{n\geq 1}\to t_0$ in  $\R$. We claim that for each $I=[q_1,q_2]$, $q_1,q_2\in\Q$, and each integer $j\geq 1$,
\begin{equation*}
	\lim_{n\to\infty}\sup_{y(\cdot)\in \K_j^I }\left|\,\int_I d(\nu_n{\cdot}{t_n})_{y(s)}-\int_I d(\nu_0{\cdot}{t_0})_{y(s)}\,\right|=0\,,
\end{equation*}
where $\K_j^I$ is the set in \eqref{eq:KjI}.

Let us fix $\ep>0$.
Since $E$ has equicontinuous $m$-bounds,  by Proposition~\ref{prop:EclsE} the closure $\overline E$ also has equicontinuous $m$-bounds, and  there is a $\delta>0$ such that  for every $\mu\in\overline E$ we have $m_j^\mu[s_1,s_2]<\ep/6$   whenever  $0\le s_2-s_1<\delta\,$. Now, take $p_1$, $p_2\in \Q$ such that $0<p_1-q_1-t_0<\delta$ and $0<q_2+t_0-p_2<\delta$. Since $t_n\to t_0$ as $n\uparrow\infty$ we can find an integer $n_0$ such that for each $n\geq n_0$,
\begin{equation}\label{equicont}
 [p_1,p_2 ]\subset  [q_1+t_n, q_2+t_n ]\,,\;\;  p_1-q_1-t_n<\delta \;\text{ and }\;  q_2+t_n-p_2<\delta\,.
\end{equation}
From Lemma~\ref{lem:integralnut} we can write
\begin{multline*}
D_n:=\sup_{y(\cdot)\in \K_j^I } \left| \int_I d(\nu_n{\cdot}{t_n})_{y(s)}- \int_ I d(\nu_0{\cdot}{t_0})_{y(s)} \right|\\
  =\sup_{y(\cdot)\in \K_j^I } \left| \int_{q_1+t_n}^{q_2+t_n} d(\nu_n)_{y(s-t_n)}- \int_{q_1+t_0}^{q_2+t_0} d(\nu_0)_{y(s-t_0)} \right|\,.
\end{multline*}
Then, decomposing the integrals using \eqref{equicont},  applying Lemma~\ref{CotasInt} and the equicontinuity of the $m$-bounds we get, for all $n\geq n_0$,
\begin{align*}
D_n &\leq  \sup_{y(\cdot)\in \K_j^I }\left| \int_{p_1}^{p_2}\left[ d(\nu_n)_{y(s-t_n)}- d(\nu_0)_{y(s-t_0)}\right]  \right| +\frac{4\,\ep}{6}\\
 &  \leq \sup_{y(\cdot)\in \K_j^I } \left| \int_{p_1}^{p_2}\left[ d(\nu_n)_{y(s-t_n)}- d(\nu_0)_{y(s-t_n)}\right]  \right|\\
 &\quad\; +\sup_{y(\cdot)\in \K_j^I } \left| \int_{p_1}^{p_2}\left[ d(\nu_0)_{y(s-t_n)}- d(\nu_0)_{y(s-t_0)}\right] \right| + \frac{2\,\ep}{3} =:P_n+ R_n +\frac{2\,\ep}{3}\,.
\end{align*}
Next, we take an interval $J$ with rational endpoints such that $I\cup [p_1,p_2]\subset J$ and, up to a suitable extension of $y(\cdot)\in \K_j^I$ by constants to $J$, the functions $y_n(\cdot)=y(\cdot-t_n)$ belong to $\K_j^J$. Therefore, since  $\{\nu_n\}_{n\geq 1}$ converges to $\nu_0$ in $\sigma_\Theta$, Lemma~\ref{lem:conv-subint} yields
\[
\lim_{n\to\infty} P_n\leq \lim_{n\to\infty} \sup_{y(\cdot)\in \K_j^J} \left|\,\int_{p_1}^{p_2} d(\nu_n)_{y(s)}-\int_{p_1}^{p_2} d(\nu_0)_{y(s)}\,\right|=0\,.
\]
Finally, Lemma~\ref{lem:xntox0} implies that $\lim_{n\to\infty} R_n=0$, which finishes the proof.
\end{proof}
\begin{rmk}\label{rm:top D}
For a set  $E\subset\PM$ with equicontinuous $m$-bounds, the continuity of the time translation map for the topology $(\PM,\sigma_D)$ can be easily deduced from the proof of Theorem~\ref{conti-time-trns}.
\end{rmk}

We conclude this section with the definition of the hull of a parametric b-measure and the construction of a continuous flow on the hull by translation.
\begin{defn}\label{def:hull}
Let $\T$ be one of the topologies $\{\sigma_\Theta, \sigma_D\}$ defined before and let $\nu\in \PM$ be a parametric b-measure.
The \emph{hull of $\nu$ with respect to $(\PM,\T)$} is defined as the topological subspace
\begin{equation*}
\mathrm{Hull}_{(\PM,\T)}(\nu)=\big(\mathrm{cls}_{(\PM,\T)}\{\nu{\cdot}t\mid t\in\R\} ,\, \T\big) \,.
\end{equation*}
\end{defn}
As a corollary of Proposition \ref{prop:EclsE}, one has the following result. We say that $\nu\in \PM$  {\it has bounded $($resp., equicontinuous$)$ $m$-bounds}  if the family $\{\nu{\cdot}t\mid t\in\R\}$ has bounded (resp., equicontinuous) $m$-bounds, and analogously with the $l$-bounds.
\begin{cor}
Let $\T$ be one of the topologies $\{\sigma_\Theta, \sigma_D\}$ and  $\nu\in \PM$.
\begin{itemize}
\item[(i)] If $\nu\in \PM$ has bounded $m$-bounds $($resp., bounded  $l$-bounds$)$, then its hull $\mathrm{Hull}_{(\PM,\T)}(\nu)$  has bounded $m$-bounds $($resp.,  bounded $l$-bounds$)$.
\item[(ii)] If $\nu\in \PM$ has equicontinuous $m$-bounds $($resp., equicontinuous  $l$-bounds$)$, then its hull $\mathrm{Hull}_{(\PM,\T)}(\nu)$  has equicontinuous $m$-bounds $($resp.,  equicontinuous $l$-bounds$)$.
\end{itemize}
\end{cor}
The final result is a corollary of Theorem~\ref{conti-time-trns}, Remark~\ref{rm:top D} and Theorem~\ref{th-compacidad}. Note that the flow condition $\varphi_0=\text{Id}$ and $\varphi_{t+s}=\varphi_t\circ\varphi_s$ for each $t, s\in\R$ ($\varphi_t(\nu)=\varphi(t,\nu)$ and $\circ$ is the composition operator) for the translation map is direct.

\begin{cor}\label{cor:hull}
Let $\nu_0\in\PM$ have equicontinuous $m$-bounds and let $\T$ be one of the topologies $\{\sigma_\Theta, \sigma_D\}$. Then, the map
\begin{equation*}
\varphi\colon\R\times \mathrm{Hull}_{(\PM,\T)}(\nu_0)\to \mathrm{Hull}_{(\PM,\T)}(\nu_0)\, ,\quad (t,\nu)\mapsto\varphi(t,\nu)=\nu{\cdot}t\,
\end{equation*}
defines a continuous flow on  $\mathrm{Hull}_{(\PM,\T)}(\nu_0)$. If in addition $\nu_0$ has bounded $l$\nbd-bounds, the hull is compact.
\end{cor}
\section{The skew-product flow induced by a class of generalized ODEs given by $N$-dim parametric b-measures on $\R^N$}\label{sec generalized ODEs}
This section is devoted to the construction of a continuous skew-product flow induced by the solutions of the equations in the hull of a generalized ODE $y'(t)=\bar \nu_{y(t)}$  given by an $N$-dim parametric b-measure $\bar\nu=(\nu^1,\ldots,\nu^N)$ on $\R^N$ with common $m$-bounds $\{m_j\}_{j\geq 1}\subset \M_c^+$ and $l$-bounds $\{l_j\}_{j\geq 1}\subset \M_c^+$ for  $1\leq i\leq N$:
\begin{itemize}
\item[(i)] $|\nu^i_y|(A)\leq m_j(A)$ for each subset $A\in \B$ and each  $y\in\R^N$ with $|y|\leq j$;
\item[(ii)] $|\nu^i_{y_1}-\nu^i_{y_2}|(A)\leq l_j(A) \,|y_1-y_2|$ for each subset $A\in\B$ and $y_1, y_2\in\R^N$ with $|y_1|, |y_2|\leq j$.
\end{itemize}
Note that, in this section, the family of moduli of continuity $\{\omega_j = {\rm id}\}_{j\geq 1}$ is necessarily related to the usual Lipschitz condition. Let us begin by defining the concept of a solution of this kind of equation.
\begin{defn}
Let $\bar\nu$ be an $N$-dim parametric b-measure on $\R^N$ with components $\bar\nu = (\nu^1,\ldots,\nu^N)$ ($\nu^i\in\PM$ for $1\leq i\leq N$) and let $t_0\in\R$ and $y_0\in\R^N$. A continuous function $y\colon[a,b] \to\R^N$, with $t_0\in [a,b]$, is said to be a \emph{solution of the Cauchy problem}
\begin{equation}\label{CauchyProblem}
\begin{cases}
y'(t)=\bar \nu_{y(t)}\\
y(t_0)=y_0
\end{cases}
\end{equation}
in $[a,b]$, if it is a solution of the integral equation  $y(t)=y_0+\int_{t_0}^td\bar\nu_{y(s)}$ for each $t\in[a,b]$, that is,
$y_i(t)=y_{0,i}+\int_{t_0}^td \nu^i_{y(s)}$ for all $t\in[a,b]$ and for each $i=1,\ldots,N$.
\end{defn}

\begin{rmk}\label{rmk:BV}
Note that, by the existence of continuous $m$-bounds and Lemma~\ref{CotasInt}, a solution is always a function of bounded variation on $[a,b]$.
\end{rmk}

\begin{exmpl}\label{ex:cantor}
A nontrivial academic example is given by the famous Cantor function $y(t)$, $t\in [0,1]$, which is the solution of the scalar generalized ODE $y'=\mu$ with initial condition $y(0)=0$, for the singular continuous Borel measure  $\mu$ on $[0,1]$ determined by the formula $\mu[0,t]=y(t)$, $t\in [0,1]$ (see, e.g., Rudin~\cite{book:Rudin}). The graph of the Cantor function, known as the Devil's staircase, is a continuous function of bounded variation which, however, is not absolutely continuous.  Recall that, by definition, the solutions of Carath\'{e}odory ODEs are within the class of absolutely continuous functions.
\end{exmpl}
Standard techniques of ordinary differential equations allow us to prove the theorem on local existence and uniqueness of solutions.
\begin{thm}\label{teor:exist uniqu sol}
Given an $N$-dim parametric b-measure $\bar\nu$ on $\R^N$ with $m$-bounds $\{m_j\}_{j\geq 1}\subset \M_c^+$ and $l$-bounds $\{l_j\}_{j\geq 1}\subset \M_c^+$ and $(t_0,y_0)\in\R\times\R^N$, there is a $\delta >0 $ such that the Cauchy problem~\eqref{CauchyProblem} has a unique solution in $[t_0-\delta,t_0+\delta]$.
\end{thm}
\begin{proof}
We consider the case $t\geq t_0$, the other case being analogous. Since the $m$-bounds and the $l$-bounds belong to $\mathcal M_c^+$, we can choose a $j\geq 2\,|y_0|\sqrt{N}$ and a $\delta>0$ such that $ m_j [t_0,t_0+\delta]\leq j/(2\sqrt{N})$ and $l_j [t_0,t_0+\delta]\leq 1/(2\sqrt{N})$. Let us consider the closed ball in $C([t_0,t_0+\delta],\R^N)$ of radius $j$,
\[
Y_j=\{y(\cdot)\in C([t_0,t_0+\delta],\R^N)\mid |y(t)|\leq j \; \text{ for all } t\in[t_0,t_0+\delta]\}\,,
\]
which is a normed space for the sup-norm. Next, for each $y(\cdot)\in Y_j$ we define $Ty$~by
\[
Ty(t)=y_0+\int_{t_0}^td\bar\nu_{y(s)}\,,\;\; t\in[t_0,t_0+ \delta]\,,
\]
and we claim that $T$ is a contractive map from $Y_j$ to $Y_j$. First notice that $Ty$ belongs to $C([t_0,t_0+\delta],\R^N)$ for each $y(\cdot)\in Y_j$ because, for each component $(Ty)_i$, from Lemma~\ref{CotasInt} applied to the component $\nu^i$, we can bound
\[
\left | (Ty)_i(t_2)-(Ty)_i(t_1)\right |=\left |\int_{t_0}^{t_2} d\nu^i_{y(s)}-\int_{t_0}^{t_1} d\nu^i_{y(s)}\right |=\left|\int_{t_1}^{t_2} d\nu^i_{y(s)}\right|\leq m_j[t_1,t_2]
\]
for each $t_0\leq t_1<t_2\leq t_0+\delta$, and $m_j\in\mathcal M_c^+$. Moreover, for every $t_0\leq t\leq t_0+\delta$,
\[
|(Ty)_i(t)|\leq |y_{0,i}|+ \left|\int_{t_0}^t d\nu^i_{y(s)}\right|\leq |y_{0,i}| + m_j[t_0,t_0+\delta]\leq \frac{j}{2\sqrt{N}}+ \frac{j}{2\sqrt{N}}=\frac{j}{\sqrt{N}}\,.
\]
From this, it follows that  $Ty\in Y_j$ and $T$ is well-defined.
If  $y(\cdot), z(\cdot)\in Y_j$,  again from Lemma~\ref{CotasInt} we deduce that for each component $1\leq i\leq N$,
\[
\left | (Ty)_i(t)-(Tz)_i(t)\right |=\left |\int_{t_0}^t d\nu^i_{y(s)}-\int_{t_0}^t d\nu^i_{z(s)}\right |\leq l_j[t_0,t_0+\delta]\,\|y-z\|_{[t_0,t_0+\delta]}\,,
\]
for each $t\in[t_0,t_0+\delta]$. Since $l_j[t_0,t_0+\delta]\leq 1/(2\sqrt{N})$, we obtain  $\|{Ty - Tz}\|_{[t_0,t_0+\delta]} \leq (1/2) \|y-z\|_{[t_0,t_0+\delta]}$. Thus, T is contractive and has a unique fixed point in $Y_j$ which is a solution of the Cauchy problem~\eqref{CauchyProblem}  in $[t_0,t_0+\delta]$. \par
To prove the uniqueness in $C([t_0,t_0+\delta],\R^N)$, if there were two solutions $y_1(\cdot)$ and $y_2(\cdot)$ of~\eqref{CauchyProblem}, define
$t_1=\sup\{t\in[t_0,t_0+\delta] \mid y_1(s)=y_2(s)  \text{ for  } t_0\leq s\leq t \}$ and let $\wt y_0=y_1(t_1)=y_2(t_1)$.
If  $t_1<t_0+\delta$, by the previous reasoning, the problem $y'(t)=\bar\nu_{y(t)}$, $y(t_1)=\wt y_0$ has a unique solution in a set
\[\wt Y=\{y(\cdot)\in C([t_1,t_1+\ep ],\R^N)\mid |y(t)|\leq k \; \text{ for all } t\in[t_1,t_1+\ep]\}\]
for large enough $k$ and small enough $\ep>0$, contradicting the definition of $t_1$.
\end{proof}
Concerning the maximal interval of definition, as in the case of Carath\'{e}odory ODEs, it is not hard to check the following result whose proof is omitted.
\begin{thm}\label{maxInt}
Given  an $N$-dim parametric b-measure $\bar\nu$ on $\R^N$  and $(t_0,y_0)\in\R\times\R^N$, there is a maximal interval $I_{\bar\nu,t_0,y_0}=(a_{\bar\nu,t_0,y_0},b_{\bar\nu,t_0,y_0})$, with $t_0\in I_{\bar\nu,t_0,y_0}$, and a continuous function  $y(t)$ defined on $I_{\bar\nu,t_0,y_0}$ which is the unique solution of \eqref{CauchyProblem}. In particular, if $a_{\bar\nu,t_0,y_0}>-\infty$ $($resp.,~$b_{\bar\nu,t_0,y_0}< \infty)$, then $|y(t)|\to \infty$ as $t\to a_{\bar\nu,t_0,y_0}^+$ $($resp.,~as $t\to b_{\bar\nu,t_0,y_0}^-)$.
\end{thm}
\begin{rmk}\label{initdata0}
    Given  an $N$-dim parametric b-measure $\bar\nu$ on $\R^N$ and $y_0\in\R^N$, we will denote by $y(\cdot,\bar\nu,y_0)$ the solution of $y'(t)=\bar\nu_{y(t)}$ with initial data $y(0)=y_0$, defined on the maximal interval of existence $I_{\bar\nu,y_0}$.
\end{rmk}
It is easy to check that the set  $\Theta$ built in the next definition, associated with a set $E$ of parametric b-measures with equicontinuous $m$-bounds, is well-defined and determines a suitable set of moduli of continuity, according to Definition~\ref{def:modcont}.
\begin{defn}\label{defimoduli}
 Let $E\subset \PM$ be a set of parametric b-measures with equicontinuous $m$-bounds and  assume  that the $m$-bounds $m_j^\nu$
  on the closed ball $B_j$ satisfy $m_j^\nu\leq m_i^\nu$
 whenever $j\leq i$ and $\nu\in E$.
 The set of functions  $\Theta=\{ \theta_j^I\colon\R^+\to\R^+ \mid I=[q_1,q_2],\,q_1,\, q_2\in \Q, \,j\in\N\}$, where
 \[
 \theta_j^I(s):=\sup_{t\in I,\, \nu\in E}\int_t^{t+s} dm_j^{\nu}\,,\quad s\geq 0\,,
 \]
 defines a suitable set of moduli of continuity associated with $E$.
\end{defn}
Next, we state the theorem of continuity of the solutions of the Cauchy problem~\eqref{CauchyProblem} with respect to initial data and the variation in $\bar\nu$. We omit the proof because it follows by adapting the steps of~\cite[Theorem~3.8~(i)]{paper:LNO2} to this case.
\begin{thm}\label{teor contin flujo}
Let $E_i\subset \PM$, $i=1,\ldots,N$ have equicontinuous $m$-bounds and consider $\Theta_i$ the suitable set of moduli of continuity built in {\rm Definition~\ref{defimoduli}} for each $i=1,\ldots,N$. Then, if  a sequence $\{\bar\nu_n\}_{n\geq 1}=\{(\nu^1_n,\ldots,\nu^N_n)\}_{n\geq 1}\subset E_1\times\ldots\times E_N$ converges to $\bar\nu\in \cls_{(\PM,\sigma_{\Theta_1})}(E_1)\times \ldots \times  \cls_{(\PM,\sigma_{\Theta_N})}(E_N) $   and $\{y_0^n\}_{n\geq 1}\subset \R^N$ converges to $y_0\in\R^N$, it holds that
\[
y(\cdot,\bar\nu_n,y_0^n) \to y(\cdot,\bar\nu,y_0) \text{ as } n\to\infty
\]
uniformly in every  $[a,b]$ contained in the maximal interval of existence $I_{\bar\nu,y_0}$.
\end{thm}
\begin{rmk}\label{rm:top D 2}
Under the hypothesis of bounded $l$-bounds,  the topologies $\sigma_{\Theta}$ and $\sigma_D$ are equivalent (see Theorem~\ref{thm:equivalenceWeakTopologies}). Therefore, with equicontinuous $m$-bounds and bounded $l$-bounds for each $E_i$ ($i=1,\ldots,N$), the previous result on continuous dependence of solutions is also true for the topology $\sigma_D$ in each component.
\end{rmk}
As in Section~\ref{sec hull}, given an $N$-dim parametric b-measure $\bar\nu_0$ on $\R^N$ with $m$-bounds $\{m_j\}_{j\geq 1}$ and $l$-bounds $\{l_j\}_{j\geq 1}$,
we say that {\it $\bar\nu_0$ has equicontinuous $m$-bounds}, if the family of translated measures
$\{m_j{\cdot}t\mid t\in\R\}$ is equicontinuous for every $j\geq 1$. Associated with this family of $m$-bounds, we consider a suitable set of moduli of continuity $\Theta$ as in {\rm Definition~\ref{defimoduli}} and the space ($\PM,\sigma_\Theta)$. Then, we endow $\PM^{\!\!\!N} $ with the product topology, $\bar\sigma_\Theta$. The hull of $\bar\nu_0$, denoted  by $\mathrm{Hull}(\bar\nu_0)$, is defined as
\begin{equation}\label{eq:hull theta}
\mathrm{Hull}(\bar\nu_0):=\big(\mathrm{cls}_{(\PM^{\!\!\!N},\bar\sigma_\Theta)}\{\bar\nu_0{\cdot}t\mid t\in\R\} ,\, \bar\sigma_\Theta\big) \,,
\end{equation}
where $\bar\nu_0{\cdot}t$ is defined componentwise, $(\nu_0^1{\cdot}t,\ldots,\nu_0^N\!{\cdot}t)$. Analogously, we say that {\it $\bar\nu_0$ has bounded $l$-bounds}, if the family of translated measures
$\{l_j{\cdot}t\mid t\in\R\}$ is bounded for every $j\geq 1$. It is easy to check that the crucial result on compactness Theorem~\ref{th-compacidad}, as well as the results in Section~\ref{sec hull}, stated for $1$-dim parametric b-measures on $\R^N$, also hold for $N$-dim parametric b-measures on $\R^N$.

Then, consider $\bar\nu_0\in\PM^{\!\!\!N}$ with equicontinuous $m$-bounds, and the family of generalized ODEs $y'(t) =\bar\nu_{y(t)}$,  $\bar\nu\in\mathrm{Hull}(\bar\nu_0)$. With the notation introduced in Remark~\ref{initdata0}, and denoting by $\U_1$   the subset of $\R\times\mathrm{Hull}(\bar\nu_0)\times\R^N$ given by
\begin{equation}\label{eq:U_1}
\U_1=\bigcup_{\substack{\bar\nu\in\mathrm{Hull}(\bar\nu_0)\\y\in\R^N}}
\{(t,\bar\nu,y)\mid t\in I_{\bar\nu,y}\}\,,
\end{equation}
 we can state the main theorem on the construction of the skew-product flow induced by the solutions of the family of generalized equations over the hull.

\begin{thm}\label{teor:conti sp}
Let $\bar\nu_0\in\PM^{\!\!\!N}$  have equicontinuous $m$-bounds $\{m_j\}_{j\geq 1}$, consider the suitable set of moduli of continuity $\Theta$ in {\rm Definition~\ref{defimoduli}} associated with the family of measures $\{m_j{\cdot}t\mid t\in\R\}$, and build the hull of $\bar\nu_0$ as in \eqref{eq:hull theta}.
 Then, the set $\,\U_1$ in \eqref{eq:U_1} is open in $\R\times\mathrm{Hull}(\bar\nu_0)\times\R^N$ and the map
\begin{equation*}
\begin{split}
\Pi\colon  \ \U_1\subset \R\times\mathrm{Hull}(\bar\nu_0)\times\R^N\  &\to\ \ \mathrm{Hull}(\bar\nu_0)\times\R^N\\
\ (t,\bar\nu,y_0)\quad \qquad \qquad &\mapsto \; (\bar\nu{\cdot}t, y(t,\bar\nu,y_0))
\end{split}
\end{equation*}
defines a local continuous skew-product flow on $\mathrm{Hull}(\bar\nu_0)\times\R^N$. If in addition $\bar\nu_0$ has bounded $l$-bounds, then the hull can be equivalently built for the product topology $\bar\sigma_D$ on $\PM^{\!\!\!N}$ for $(\PM,\sigma_D)$ in each component, and the hull is compact.
\end{thm}
\begin{proof}
First of all, $\Pi$ is a skew-product flow (or nonautonomous dynamical system) because we have a continuous flow on the base after Corollary~\ref{cor:hull} and, using Lemma~\ref{lem:integralnut}, the fiber components satisfy the cocycle property, subject to existence:
\begin{equation*}
  y(t+s,\bar\nu,y_0)=y(t,\bar\nu{\cdot}s,y(s,\bar\nu,y_0))\,,\quad t,s\in\R\,,\;  (\bar\nu,y_0)\in \mathrm{Hull}(\bar\nu_0)\times\R^N.
\end{equation*}
The continuity and the remaining assertions are a consequence of Theorem~\ref{thm:equivalenceWeakTopologies},  Theorems~\ref{maxInt} and~\ref{teor contin flujo}, and Remark~\ref{rm:top D 2}. The proof is finished.
\end{proof}
\subsection{Application to precompact families of Carath\'{e}odory ODEs}\label{subsec appl Carath}
The previous skew-product framework allows us to apply topological dynamics techniques to the study of nonautonomous ODEs $y'=g(y,t)$ with Carath\'{e}odory conditions, as introduced at the beginning of Section~\ref{sec Motivation}: namely, recall conditions (i) and (ii) on the map $g$, and the properties of equicontinuous $m$-bounds and bounded $l$-bounds in (i)$'$ and (ii)$'$, to which we will refer in what follows. In particular, note that the maps $\wt m_K, \wt l_K \in L^1_{loc}$ are also $m$-bounds/$l$-bounds of each component of $g=(g_1,\ldots,g_N)$.

The process is as follows. We associate an $N$-dim parametric b-measure on $\R^N$ with the map $g(y,t)$. Consider $\bar\nu\colon\R^N\to \M_c(\B)^N$, $y\mapsto (\nu_y^1,\ldots, \nu_y^N)$ given by
\[
\nu_y^i(A):= \int_A g_i(y,t)\,dt\quad\hbox{for each } A\in \mathcal{B} \;\, (1\leq i\leq N)\,,
\]
which is well-defined by (i) and the fact that $\wt m_K\in L^1_{loc}$. Notice that, using a differential notation, we can simply write
\begin{equation}\label{eq:diff not}
d\bar \nu_y = g(y,t)\,dt \quad\hbox{for each }\, y\in \R^N.
\end{equation}
Let us check that $\bar\nu$ is an $N$-dim parametric b-measure on $\R^N$ with equicontinuous $m$-bounds and bounded $l$-bounds.
For each integer $j\geq 1$, let $m_j$ and $l_j$ be the absolutely continuous regular Borel measures with density functions $\wt m_{B_j} (t)$ and $\wt l_{B_j} (t)$, respectively. Then, $m_j, l_j\in \M_c^+$ and it is easy to check that
\begin{itemize}
\item[(1)] $|\nu^i_y|(A)\leq m_j(A)$ for each subset $A\in \B$ and each  $y\in\R^N$ with $|y|\leq j$,
\item[(2)] $|\nu^i_{y_1}-\nu^i_{y_2}|(A)\leq l_j(A) \,|y_1-y_2|$ for each subset $A\in\B$ and $y_1, y_2\in\R^N$ with $|y_1|, |y_2|\leq j$,
\end{itemize}
and this holds for each $1\leq i\leq N$. That is, $m_j, l_j\in \M_c^+$ are, respectively, an $m$-bound and an $l$-bound of $\nu^i$ on the ball $B_j$ according to Definition~\ref{defi:parameasure}. It is straightforward that the desired properties of equicontinuity of the set of translated measures $\{m_j{\cdot}t\mid t\in\R\}$ and boundedness of $\{l_j{\cdot}t\mid t\in\R\}$ for each $j\geq 1$ are inherited from (i)$'$ and (ii)$'$ cited
above.  Furthermore, with bounded $l$-bounds, we can build the hull of $\bar\nu$ using the topology $\sigma_D$ on $\PM$ (see Theorem~\ref{thm:equivalenceWeakTopologies}), which is easier to handle than the topology $\sigma_\Theta$. A fundamental result is the following.
\begin{prop}\label{prop:same solutions}
Suppose that $g\colon\R^N\times\R\to \R^N$ is measurable in $t$, continuous in $y$ for almost every $t\in\R$, and for every compact set $K\subset \R^N$ there are $m$-bounds, $l$-bounds, $\wt m_K, \wt l_K\in L^1_{loc}$, respectively. Then, $y(t)$ is a solution of the Carath\'{e}odory ODE $y'=g(y,t)$  if and only if it is a solution of the associated generalized ODE $y'(t)=\bar\nu_{y(t)}$, for the parametric b-measure $\bar\nu$ in \eqref{eq:diff not}.
\end{prop}
\begin{proof}
It suffices to prove that, if $y\colon[a,b]\to \R^N$ is continuous, then the integral of $\bar\nu$ along $y\colon[a,b]\to \R^N$, $\int_a^b d\bar\nu_{y(s)}$  coincides with the Lebesgue integral $\int_a^b g(y(s),s)\,ds$. To see it,  for each $n\geq 1$ take $\delta_n=1/n$  and a tagged $\delta_n$-partition
\[
\Delta^{\delta_n}= \{ (\tau_i^n,[t_i^n,t_{i+1}^n]) \mid i=1,\ldots, k_n-1\}\,,
\]
and consider the map $g_n\colon[a,b]\to \R^N$, $g_n(s):= \sum_{i=1}^{k_n-1} \chi_{[t_i^n,t_{i+1}^n)}(s)\,g(y(\tau_i^n),s)$ for $s\in [a,b)$ and $g_n(b):=g(y(\tau_{k_n-1}^n),b)$. By construction, $g_n$ is Lebesgue measurable. Since $g(\cdot,s)$ is continuous for almost every $s\in[a,b]$, and $y\colon[a,b]\to \R^N$ is continuous, it easily follows that for almost every $s\in [a,b]$ there exists the pointwise limit $\lim_{n\to\infty} g_n(s) = g(y(s),s)$. Furthermore, if $K\subset \R^N$ is a compact set such that $y([a,b])\subset K$, then for each $n\geq 1$, $|g_n(s)|\leq \wt m_K(s)$ for a.e.~$s\in[a,b]$, and $\wt m_K\in L^1([a,b])$.  Thus, Lebesgue's dominated convergence theorem implies $\int_a^b g(y(s),s)\,ds = \lim_{n\to\infty} \int_a^b g_n(s)\,ds$. Since
\[
\int_a^b g_n(s)\,ds=\sum_{i=1}^{k_n-1} \int_{t_i^n}^{t_{i+1}^n} g(y(\tau_i^n),s)\,ds= \sum_{i=1}^{k_n-1} \int_{t_i^n}^{t_{i+1}^n} d\bar\nu_{y(\tau_i^n)} = S_{\bar\nu}^{y(\cdot)}(\Delta^{\delta_n})
\]
for the Riemann sums introduced in Section~\ref{sec Integ}, and $\delta_n\to 0$, the limit is equal to $\int_a^b d\bar\nu_{y(s)}$. The proof is finished.
\end{proof}
In summary, the ODE  $y'=g(y,t)$ considered in this section can be rewritten as the generalized ODE given by the $N$-dim parametric b-measure $\bar\nu$ on $\R^N$ in~\eqref{eq:diff not}, that is, $y'(t)=\bar\nu_{y(t)}$, and it can receive a nonautonomous dynamical treatment under the skew-product formalism explained in this section. We highlight two facts: first, the flow on the hull is compact (see Theorem~\ref{teor:conti sp}), and second,  all the equations over the hull keep the properties of equicontinuous $m$-bounds and bounded $l$-bounds (see Proposition~\ref{prop:EclsE}) and their solutions are all continuous functions of bounded variation  (see Remark~\ref{rmk:BV}). This avoids taking a larger set of generalized ODEs where solutions lacking these properties may appear. 
\subsection{Numerical examples}\label{subsec:numerics}
We aim to provide some examples of how solutions of the type of differential equation described in this work may look in practice.

As a basis for our examples, we use an approximation of the Cantor function on [0,1], whose graph is the well-known Devil's staircase. We approximate the Cantor function by integrating the maps $f_n(t)$ built as follows. We consider the $n$-th iteration of the middle-third Cantor set, which consists of $2^{n}$ disjoint intervals within $[0,1]$, $\{I_1^n,\ldots,I_{2^n}^n\}$, all of them of length $1/3^{n}$. A piecewise constant function $f_n(t)$ is then defined to be equal to a constant value $a_n$ on these intervals and zero elsewhere, where $a_n$ is chosen so that the integral of $f_n$ over $[0,1]$ is normalized to one:
\[
f_n(t)=\sum_{j=1}^{2^{n}} a_n\, \chi_{I^n_j}(t)\,.
\]
In particular, we will consider $n=12$ in the numerical simulations below. For $n=12$, we have $a_{12}\approx 129.7463$.

To introduce temporal complexity, we consider linear combinations of translated versions of $f_n(t)$. Firstly, we define
\begin{equation}\label{eq:cantor-translations}
f_n^{\textit sum}(t) = \sum_{k=1}^{100} (-1)^k f_n(t - k\delta)\,,
\end{equation}
where $\delta = 0.01$. This alternating sum introduces oscillatory behavior while preserving the fractal support structure. The maps $f_n^{\textit sum}(t)$ are Carath\'{e}odory maps with compact support. More precisely,  $f_n^{\textit sum}(t)=0$ if $t\leq 0$ or $t>2$ and $f_n^{\textit sum}(2)=f_n(1)=a_n$.  Now, $\int_0^t f_n^{\textit sum}(s)\,ds$ converges uniformly on $[0,2]$ to a linear combination of singular continuous measures $\sum_{k=1}^{100} (-1)^k \mu{\cdot}(-k\delta) [0,t]$, given by translations of the singular continuous measure $\mu$ associated with the Cantor function (see Example~\ref{ex:cantor}),  as $n\to\infty$.

\begin{figure}[h]
    \centering
    \begin{overpic}[width=\textwidth]{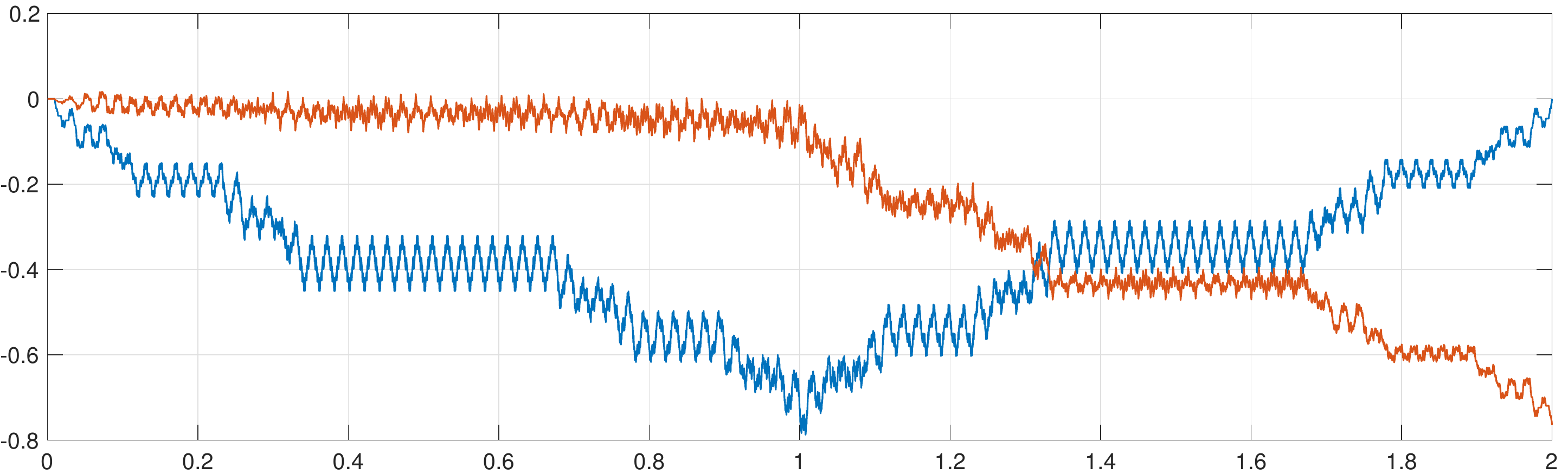}
    \put(50,-2){$t$}
    \put(4.5,26.5){$y$}
    \end{overpic}
    \vspace{-1mm}

    \begin{overpic}[width=\textwidth,trim={-0.9cm 0 0 -5mm},clip]{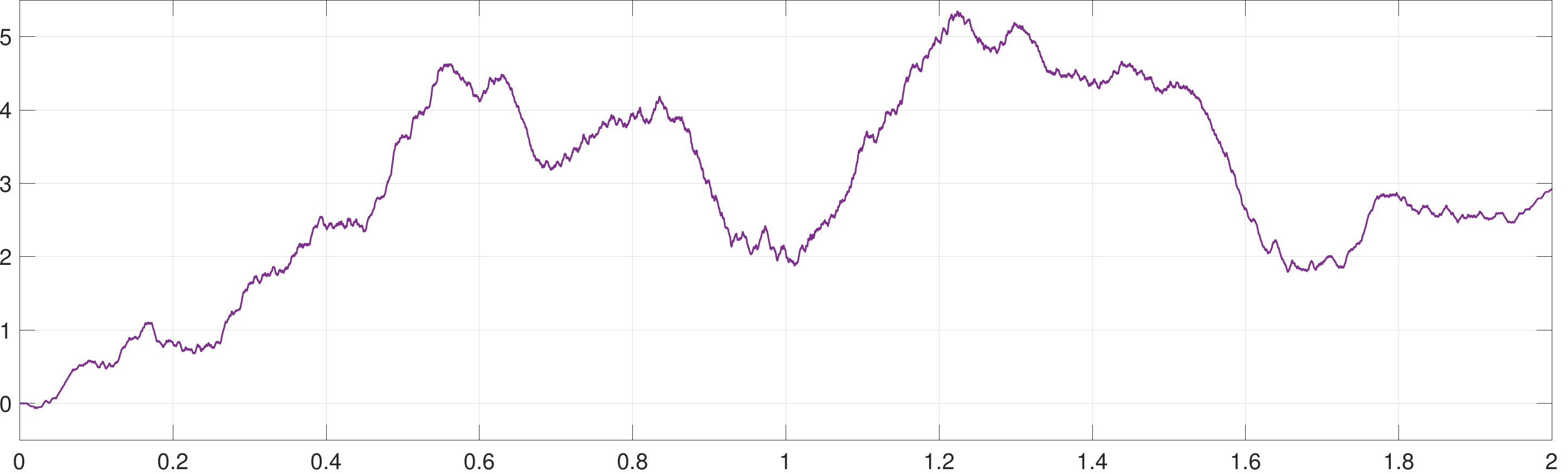}
    \put(50,-2){$t$}
    \put(4.5,26.5){$y$}
    \end{overpic}

    \caption{Upper panel: in solid blue is depicted the graph of the solution to the initial value problem $ y'=f_{12}^{\textit sum}(t)$, $y(0)=0$, where $f_{12}^{\textit sum}(t)$ is defined as in \eqref{eq:cantor-translations}. In solid orange the solution to the initial value problem with right-hand side defined via $\overline f_{12}^{\,{\textit sum}}(t)$ (cf.~\eqref{eq:cantor-translations2}). Lower panel: graph of the solution to the initial value problem $ y'=\widetilde f_{12}^{\textit sum}(t)$, $y(0)=0$ where $\widetilde f_{12}^{\,{\textit sum}}(t)= \sum_{k=1}^{100} \textrm{sign}(2\,\textrm{rand()-1})f_{12}(t - k\delta)$, that is, the coefficients of the translations of $f_{12}$ are selected randomly in $\{-1,0,1\}$.\\[1ex]}
    \label{fig:Cantor1}
\end{figure}

The solid blue line in the upper panel of  Figure~\ref{fig:Cantor1} is the graph of the function
\[
y(t) =\int_0^t f_{12}^{\textit sum}(s)\,ds\,, \quad  t\in [0,2]\,,
\]
which is the solution of the initial value problem $y'= f_{12}^{\textit sum}(t)$, $y(0)=0$. Note that the Carath\'{e}odory ODE $y'=f_{12}^{\textit sum}(t)$ is an approximation of the generalized ODE $y'= \sum_{k=1}^{100} (-1)^k \mu{\cdot}(-k\delta)$ for $t\in [0,2]$, and recall that solutions depend continuously on variations of the equation, as stated in Theorem~\ref{teor contin flujo}.
In fact, we can build a Carath\'{e}odory ODE in the class that motivates this paper (that is, such that the hull construction scapes the context of ODEs) so that the map $y(t)$ is a solution on $[0,2]$. For this purpose, we  fix a sequence of integers $\{n_j\}_{j\geq 0}$,  with $n_0=12$ for convenience, and such that  $n_j\uparrow\infty$, and consider the map
\begin{equation}\label{eq:precompact map h}
    h(t) = \sum_{j=0}^\infty f_{n_j}^{\textit sum}(t-2j)\,.
\end{equation}
For $t\in [0,2]$, $h(t) =f_{12}^{\textit sum}(t)$, since $f_{n_j}^{\textit sum}(t-2j)=0$ for every $j\geq 1$. Thus, $y(t)$ satisfies $y'=h(t)$ for $t\in [0,2]$. Futhermore, the equation $y'=h(t)$ has equicontinuous $m$-bounds and in its hull there are singular continuous measures: note that the translates of $h$,  $h(t+2k)$ ($k\geq 1$) are such that the absolutely continuous measures with densities $h(t+2k)$ restricted to $[0,2]$ coincide with the measures with densities $f_{n_k}^{\textit sum}(t)$, which converge to a singular continuous measure on  $[0,2]$ as $k\to\infty$, by construction.

Furthermore, for integration on a compact interval, the choice of coefficients $(-1)^k$ in the linear combination \eqref{eq:cantor-translations} is absolutely arbitrary. For example, the orange solid line in Figure \ref{fig:Cantor1} is the graph of the solution to the initial value problem $y' = \overline f_{12}^{\,{\textit sum}}(s)\,ds$, with $y(0)=0$ and
\begin{equation}\label{eq:cantor-translations2}
\overline f_n^{\,{\textit sum}}(t) = \sum_{k=1}^{100} \sin\left(k\sqrt{11}\right) f_n(t - k\delta)\,,
\end{equation}
while the lower panel depicts the solution to the analogous initial value problem when the coefficient is randomly chosen in $\{-1,0,1\}$.

The function $f_{12}^{\textit sum}(t)$ in \eqref{eq:cantor-translations} is also used as a time-dependent forcing term in the nonautonomous Riccati differential equation
\[
y' = -y^2 + 2 + f_{12}^{\textit sum}(t)\,.
\]
A complete dynamical description of nonautonomous Riccati equations with piecewise uniformly continuous coefficients can be found in~\cite{paper:LNuO}. In particular, if this type of equation has uniformly separated bounded solutions, then it has exactly two hyperbolic solutions, the upper one stable and the lower one unstable, which delimit the set of bounded solutions of the equation respectively from above and from below.
\begin{figure}[h]
    \centering
    \begin{overpic}[width=\textwidth]{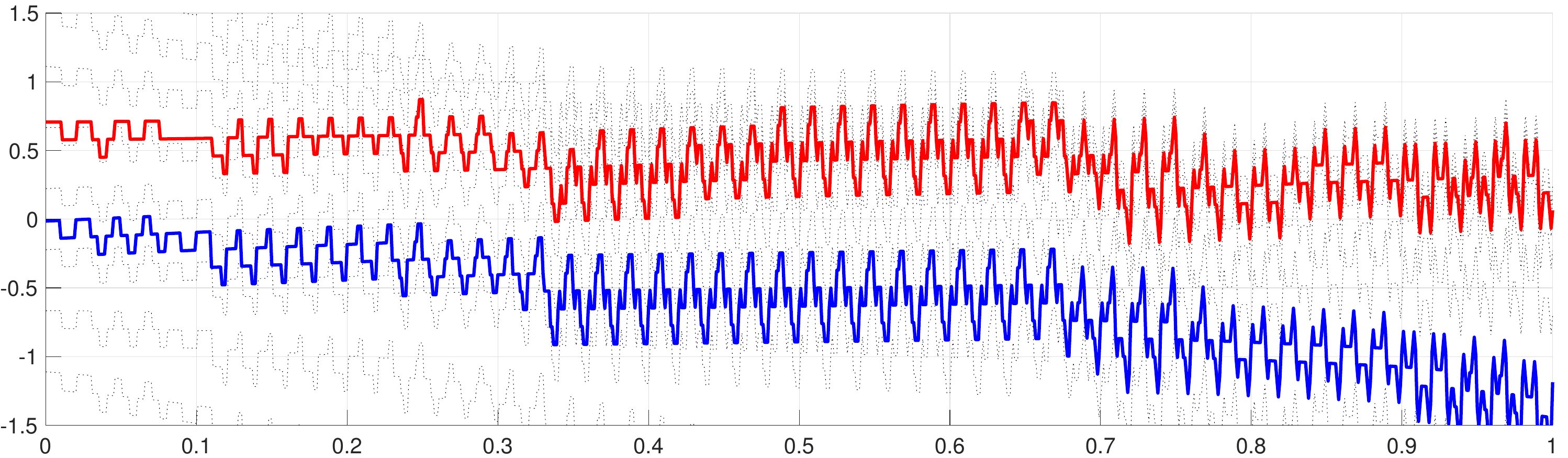}
    \put(50,-2){$t$}
    \put(4.5,26.5){$y$}
    \end{overpic}
    \caption{A selection of solutions of the differential equation $y' = -y^2+ 2 + f_{12}^{\textit sum}(t)$ are plotted in the above panel. The red curve approximates the graph of the hyperbolic attracting solution, while the blue one approximates the graph of the hyperbolic repelling solution. The black dotted curves are examples of eight additional solutions with different initial conditions at time zero.\\[1ex]}
    \label{fig:Cantor}
\end{figure}
Figure~\ref{fig:Cantor} displays numerical solutions of the Riccati equation using the Euler's method (chosen instead of a Runge-Kutta method due to the stiff nature of the differential problem at hand): eight trajectories (in black) starting at $t=0$ with initial conditions uniformly spaced in $[-2,2]$, one trajectory (in red) integrated forward from $t=-100$ with $y(-100)=2$  approximating the unique hyperbolic stable solution, and one (in blue) integrated backward from $t=100$ with $y(100)=-2$  approximating the unique hyperbolic repelling solution. For the numerical simulations, the function $h$ in \eqref{eq:precompact map h} is considered only for $t\le 2$ and it thus coincides with $f_{12}^{sum}$.

\section{Application to slow-fast systems of ODEs}\label{sec slowfast}
Recently, Longo et al.~\cite{paper:LOStracking} have obtained a dynamical interpretation of the behaviour of the fast variables of slow-fast systems of ODEs with dependence on the fast time in terms of tracking of nonautonomous attractors. In this section, we review the validity of their results under weaker regularity conditions on the vector field driving the fast motion, such as those imposed in Section~\ref{sec Motivation}, making the field of application to models in real life much wider.

Keeping the notation of the cited paper, in this type of singularly perturbed system there are two scales of time: the slow timescale $t$ and the fast timescale $\tau$; and also two scales of motion: the slow variables $x\in\R^n$ and the fast variables $y\in\R^m$. Written in slow time, $t$, the coupled pair of ODEs takes the form
\begin{equation}\label{eq:slowfastent}
\left\{\begin{array}{l} \dot{x}  = f(x,y)\,,
  \\
  \ep\,\dot{y}  = g\left(x,y,t/\ep \right),
\end{array}\right.
\end{equation}
where $\dot{x}, \dot{y}$ denote the derivatives with respect to $t$ and $\ep>0$ is a small parameter. After the change of independent variables $\tau=t/\ep$ relating the two time scales, the equations are written as
\begin{equation}\label{eq:slowfastentau}
\left\{\begin{array}{l} x'  = \ep\,f(x,y)\,,
  \\
 y' = g(x,y,\tau)\,,
\end{array}\right.
\end{equation}
where $x', y'$ denote the derivatives with respect to $\tau$.
 Note that the equations for the fast motion explicitly depend on the fast time, which makes the problem nonautonomous. Prescribed initial conditions are given:
\begin{equation}\label{eq:initial values}
\left\{\begin{array}{l} x(0)=x_0\,,
  \\
  y(0)=y_0\,.
\end{array}\right.
\end{equation}
For each fixed $\ep>0$, the solution of the initial value problem~\eqref{eq:slowfastent}-\eqref{eq:initial values} will be denoted by $(x_\ep(t),y_\ep(t))$, for $t\in[0,t_0]$, and the solution of~\eqref{eq:slowfastentau}-\eqref{eq:initial values} will be denoted by $(x_\ep(\tau),y_\ep(\tau))$, for $\tau\in[0,t_0/\ep]$. The main interest is in the behaviour of the solutions as $\ep\to 0$, and specifically in the limit behaviour of $y_\ep(t_0/\ep)$. We refer the reader to Tikhonov~\cite{paper:Tikh}, Verhulst~\cite{book:Ve}, Berglund and Gentz~\cite{book:BG}, and Kuehn~\cite{book:ku} for the theory and applications of autonomous slow-fast systems.  When dependence on the fast time is admitted, the first work is due to  Artstein~\cite{paper:ZA99}, followed by the more recent results in~\cite{paper:LOStracking} with a different dynamical approach to the behaviour of the fast motion.

The basic assumptions on the vector fields $f$ and $g$ are as follows. Note that, compared to Assumption~3.1 in~\cite{paper:LOStracking}, assumption~(i) on $f$ is the same, but assumptions~(ii)-(iii) on $g$ are much weaker than~(ii)-(iii) in~\cite{paper:LOStracking}, where $g$ is bounded and uniformly continuous in $\tau$, and uniformly Lipschitz in $y$ on compact sets.

Due to the presence of two types of variables, in this section the closed ball of $\R^N$ centered at the origin and with radius $j$ is denoted by $B_j^{(N)}$.
\begin{assu}\label{asu:1}
\begin{itemize}
  \item[]
  \item[(i)] $f(x,y)$ is continuous and $\sup_{y\in \R^m} |f(x,y)|\leq a + b\,|x|$
  for all $x\in\R^n$, for some  $a$ and $b$ positive.
  \item[(ii)] $g(x,y,\tau)$ is measurable in $\tau$, continuous in $(x,y)$ for a.e.~$\tau$, and for each integer $j\geq 1$  there is a modulus of continuity $\omega_j\colon\R^+\to \R^+$, with $\omega_j\leq \omega_{j+1}$, and two nonnegative locally integrable maps $\wt m_j, \wt l_j \in L^1_{loc}$ (respectively called an $m$-bound  and  an $l$-bound of $g$  on $B_j^{(n)}\times B_j^{(m)}$)  such that:
\begin{itemize}
    \item[(a)] for a.e.~$\tau\in\R$,  $|g(x,y,\tau)| \leq  \wt m_j (\tau)$ for all $x\in B_j^{(n)}$ and  $y\in B_j^{(m)}$;
    \item[(b)] for a.e.~$\tau\in\R$,  $|g(x_1,y_1,\tau)-g(x_2,y_2,\tau)| \leq  \wt l_j (\tau)\,[\omega_j(|x_1-x_2|)+|y_1-y_2|]$ for all $x_1,x_2\in B_j^{(n)}$ and  $y_1,y_2\in B_j^{(m)}$.
\end{itemize}
\item[(iii)] The $m$-bounds $\wt m_j$ are equicontinuous and the $l$-bounds $\wt l_j$ are bounded, that is, they respectively satisfy (i)$'$ and (ii)$'$ in Section~\ref{sec Motivation}.
\end{itemize}
\end{assu}
Since the problem for the fast variables $y'=g(x,y,\tau)$ is nonautonomous, our approach is to build a nonautonomous forward dynamical system or skew-product semiflow over a compact base flow. In this context, the base flow must take into account the closure of the time-translated problems $y' = g(x,y,\tau_0+\tau)$ $(\tau_0\in \R)$  as well as the variation of the parameter $x$. As discussed in Section~\ref{sec generalized ODEs}, to this end, we can rewrite the ODEs as generalized ODEs associated with an $m$-dim parametric b-measure on $\R^{n}\times \R^m$. More precisely, associated with the vector field $g$ we build an $m$-dim parametric b-measure $\bar\nu_g\colon\R^{n}\times \R^m\to \M_c(\B)^m$, $(x,y)\mapsto (\bar\nu_g)_y^x$, with components $(\nu_g^1,\ldots,\nu_g^m)$, written in differential form as
\[
d(\bar\nu_g)_y^x = g(x,y,\tau)\,d\tau\,, \quad x\in \R^n,\;y\in\R^m.
\]
For each integer $j\geq 1$, we take $m_j$ and $l_j$ the absolutely continuous regular Borel measures with density functions $\wt m_j (\tau)$ and $\wt l_j (\tau)$, respectively.  Then, $m_j, l_j\in \M_c^+$ and it is easy to check that for each component $1\leq i\leq m$,
\begin{itemize}
\item[(1)] $|(\nu_g^i)_y^x|(A)\leq m_j(A)$ for each subset $A\in \B$ and each $(x,y)\in B_j^{(n+m)}$,
\item[(2)] $|(\nu_g^i)_{y_1}^{x_1}-(\nu_g^i)_{y_2}^{x_2}|(A)\leq l_j(A) \,[\omega_j(|x_1-x_2|)+|y_1-y_2|]$ for each subset $A\in\B$ and $(x_1,y_1),\,(x_2,y_2)\in B_j^{(n+m)}$.
\end{itemize}
That is,  $m_j, l_j\in \M_c^+$ are, respectively, an $m$-bound and an $l$-bound of $\nu_g^i$ on the ball $B_j^{(n+m)}$, with respect to a particular family of moduli of continuity, which is Lipschitz in the $y$ variables. Therefore, the results in Sections~\ref{sectopo} and~\ref{sec hull} apply, and those in Section~\ref{sec generalized ODEs} apply to $(\bar\nu_g)^x\colon\R^m\to \M_c(\B)^m$, $y\mapsto (\bar\nu_g)_y^x$,  for each $x$ fixed.

Note that we assign different roles to each of the variables $x$ and $y$. The slow variables $x$ play the role of a vector-valued parameter, whereas the fast variables $y$ are those with respect to which we solve a generalized ODEs problem. More precisely, for each $x\in\R^n$ fixed, the equation $y' = g(x,y,\tau)$ is called a {\it layer problem\/} or a {\it layer equation\/} which is rewritten (see Proposition~\ref{prop:same solutions}) as
\[
y'(\tau) = (\bar\nu_g)_{y(\tau)}^x\,,\quad \tau\geq 0\,,
\]
and the translated problems as $y'(\tau) = (\bar\nu_g{\cdot}\tau_0)_{y(\tau)}^x$ $(\tau_0\in \R)$.
Since, according to Assumption~\ref{asu:1}~(iii), the $m$-dim parametric b-measure $\bar\nu_g$ on $\R^n\times\R^m$ has equicontinuous $m$-bounds and bounded $l$-bounds,  we have a continuous flow defined by translation on the compact hull $\mathcal{H}$ of $\bar\nu_g$ built as in \eqref{eq:hull theta}, which can equivalently be built for the product topology $\bar\sigma_D$.  We build an extended compact base flow by adding a compact component $K_0\subset \R^n$,
\begin{equation*}
 \begin{array}{ccl}
 \R\times\mathcal{H}\times{K_0}& \longrightarrow & \mathcal{H}\times{K_0}\\
 (\tau,\bar\nu,x) & \mapsto &(\bar\nu{\cdot}\tau,x)\,.
\end{array}
\end{equation*}
Then, under Assumption~\ref{asu:1}~(ii)-(iii), a continuous skew-product semiflow is defined over the former compact base flow in the following way, on an appropriate open set $\mathcal{U}$ subject to the existence of solutions:
\begin{equation}\label{eq:sk-pr}
 \begin{array}{ccl}
 \pi\colon\mathcal{U}\subseteq\R^+\times\mathcal{H}\times{K_0} \times \R^m & \longrightarrow & \mathcal{H}\times{K_0} \times \R^m \\
 (\tau,\bar\nu,x,y_0) & \mapsto &(\bar\nu{\cdot}\tau,x,y(\tau,\bar\nu,x,y_0))\,,
\end{array}
\end{equation}
where $y(\tau,\bar\nu,x,y_0)$ denotes the solution of $y'(t)=\bar\nu_{y(t)}^x$ with initial data $y(0)=y_0$. Keep in mind that for the element $\bar\nu_g\in\mathcal{H}$ we recover the solutions of the Carath\'{e}odory ODEs $y'=g(x,y,\tau)$ for $x\in K_0$ (once again, see Proposition~\ref{prop:same solutions}).

We insist that the Lipschitz-type condition in Assumption~\ref{asu:1}~(ii)-(b) with respect to the variable $y$ is crucial to apply the results in Section~\ref{sec generalized ODEs}. In particular, Theorem~\ref{teor contin flujo} and Remark~\ref{rm:top D 2} assert the continuous dependence of the solutions with respect to the variation of $\bar\nu\in \mathcal{H}$ and the initial condition $y_0$. To check the continuous dependence with respect to changes in the parameter $x$, the arguments are the same (once more, see the proof of~\cite[Theorem~3.8~(i)]{paper:LNO2}), and we omit them.
Note that $\pi$ is a nonautonomous dynamical system because, by construction, the fiber components satisfy the so-called {\it semicocycle property}, subject to existence:
\begin{equation}\label{eq:cocycle}
  y(\tau_1+ \tau_2,\bar\nu,x,y_0)=y(\tau_1,\bar\nu{\cdot}\tau_2,x,y(\tau_2,\bar\nu,x,y_0))\,,\quad \tau_1,\tau_2\geq 0\,.
\end{equation}

For each fixed $x\in\R^n$, let us denote by $\pi^x$ the {\it layer skew-product semiflow\/}:
\begin{equation}\label{eq:sk-pr x}
\begin{split}
\pi^x\colon  \mathcal{U}_x\subseteq \R^+\times\mathcal{H}\times\R^m\  &\to\ \ \mathcal{H}\times\R^m\\
\ (\tau,\bar\nu,y_0)\qquad &\mapsto \; (\bar\nu{\cdot}\tau, y(\tau,\bar\nu,x,y_0))\,.
\end{split}
\end{equation}

As in~\cite{paper:LOStracking}, the main hypothesis in qualitative behaviour of the fast variables is the uniform ultimate boundedness of the solutions of the parametric family of Carath\'{e}odory ODEs $y' = g(x,y,\tau)$ for $x$ varying in compact subsets $K_0\subset \R^n$. This requires that solutions are globally defined in the future and that there exists a constant $c>0$ such that  if initial conditions $y_0$ are taken in a ball $B_d\subset \R^m$, there exists a time $T=T(d)>0$ such that the solutions have entered the ball $B_c\subset \R^m$ after a lapse of time $T$, and this happens uniformly for $x\in K_0$ (more precisely, see Definition~3.2 in~\cite{paper:LOStracking}). We collect these additional hypotheses.
\begin{assu}\label{asu:2}
\begin{itemize}
  \item[]
  \item[(i)] For each $\bar\nu\in\mathcal{H}$, $x\in \R^n$  and $y_0\in \R^m$, the solution $y(\tau,\bar\nu,x,y_0)$ is defined for all $\tau\geq 0$.
  \item[(ii)] For each compact set $K_0$ of $\R^n$, the solutions of the parametric family $y'=g(x,y,\tau)$,  $x\in K_0$ are uniformly ultimately bounded.
\end{itemize}
\end{assu}
Once we have established the context, it is important to note that Theorem~5.3 in Artstein~\cite{paper:ZA99}, which basically states uniform convergence for the slow variables and statistical convergence for the fast variables, applies under stronger assumptions on  $g$:
the uniform integrability of the $m$-bounds is now relaxed to equicontinuity, with the technical, though challenging, inconvenience that generalized ODEs appear, but with the advantage of enlarging the field of applications.  We keep their result for the slow variables, whereas we follow the dynamical approach of Longo et al.~\cite{paper:LOStracking} for the qualitative behaviour of the fast variables, using the theory of nonautonomous attractors. This allows us to identify the dynamical objects that the fast variables approach as $\ep$ goes to $0$, in terms of the fibers of the global attractor.

Not to make the paper unnecessarily long, we refer the reader to~\cite{paper:LOStracking} and references therein (e.g., Kloeden and Rasmussen~\cite{klra} and~\cite[Section 5]{paper:LNO2} in a Carath\'{e}odory context) for the concepts of global attractor, pullback attractor and parametrically inflated pullback attractor, which will appear in the statements of the results. For the proofs of the next two results, the reader is referred to the proofs of~\cite[Proposition 3.4]{paper:LOStracking} and~\cite[Theorem~3.9]{paper:LOStracking}, respectively.
\begin{prop}\label{prop:existe atractor}
Under Assumptions~{\rm\ref{asu:1}~(ii)-(iii)}  and~{\rm\ref{asu:2}},  given a compact set $K_0\subset \R^n$, the following statements hold:
\begin{itemize}
  \item[\rm{(i)}] There is a global attractor $\A\subset \mathcal{H}\times{K_0} \times \R^m$ for the globally defined skew-product semiflow $\pi$ given in~\eqref{eq:sk-pr} which can be written as
\[
\A=\bigcup_{(\bar\nu,x)\in\mathcal{H}\times{K_0}} \{(\bar\nu,x)\}\times A_{(\bar\nu,x)}
\]
for the sections $A_{(\bar\nu,x)}=\{ y\in \R^m\mid (\bar\nu,x,y)\in \A \}$, and the nonautonomous set $\{A_{(\bar\nu,x)}\}_{(\bar\nu,x)\in \mathcal{H}\times{K_0}}$ is the pullback attractor of $\pi$.
\item[\rm{(ii)}] For each fixed $x\in K_0$, the skew-product semiflow $\pi^x$ given in~\eqref{eq:sk-pr x} has a global attractor $\A^x\subset \mathcal{H}  \times \R^m$ which can be written as
\[
\A^x=\bigcup_{ \bar\nu\in\mathcal{H}} \{\bar\nu\}\times A^x_{\bar\nu}
\]
for the sections $A^x_{\bar\nu}=\{ y\in \R^m\mid (\bar\nu,y)\in \A^x \}$, and the nonautonomous set $\{A^x_{\bar\nu}\}_{\bar\nu\in \mathcal{H}}$ is the pullback attractor of the skew-product semiflow $\pi^x$.
\end{itemize}
   \end{prop}
For $\delta>0$, the nonautonomous set $\{A_{(\bar\nu,x)}[\delta]\}_{(\bar\nu,x)\in \mathcal{H}\times K_0}$ is the parametrically inflated pullback attractor of the skew-product semiflow $\pi$ (see~\cite{klra},~\cite{paper:LOStracking}).

\begin{thm}\label{teor:forward}
Under Assumptions~{\rm\ref{asu:1}~(ii)-(iii)}  and~{\rm\ref{asu:2}},  given a compact set $K_0\subset \R^n$, let $\A$ be the global attractor given in Proposition~$\ref{prop:existe atractor}~\rm{(i)}$. Then, fixed a $\delta>0$,  for every bounded set $D\subset \R^m$,
\begin{equation*}
  \lim_{\tau\to\infty} \bigg( \sup_{\bar\nu\in \mathcal{H},\, x\in K_0} {\rm dist}\big(y(\tau,\bar\nu,x,D),A_{(\bar\nu{\cdot}\tau,x)}[\delta]\big) \bigg) =0\,,
\end{equation*}
for the Hausdorff semidistance ${\rm dist}$.
\end{thm}
In this Carath\'{e}odory context, due to the relaxed Lipschitz-type condition in Assumption~\ref{asu:1}~(ii)-(b), which depends on the $l$-bound maps, we will use the following result instead of Lemma~4.1 in~\cite{paper:LOStracking}.
\begin{lem}\label{lem:comparacion}
Let $g_1,g_2\colon\R^m\times \R\to \R^m$ be Carath\'{e}odory maps with $m$-bounds, and assume that $g_1$ also has bounded $l$-bounds $\{\wt \ell_j\}_{j\geq 1}\subset L^1_{loc}$. If $y_1,\, y_2\colon[a,b]\to \R^m$ are solutions of $y'=g_1(y,\tau)$, $y'=g_2(y,\tau)$, respectively, with $y_1(a)=y_2(a)$, $y_1([a,b]),y_2([a,b])\subset B_j^{(m)}$ for some $j\geq 1$, and for some $\sigma>0$,
\begin{equation}\label{eq:aux 4}
|g_1(y_2(\tau),\tau)-g_2(y_2(\tau),\tau)|\leq  \wt \ell_j(\tau)\,\sigma \quad \hbox{for almost every }\;\tau\in [a,b]\,,
\end{equation}
 then, for a positive constant $c_j$ such that $\int_s^{s+(b-a)} \wt \ell_j(\tau)\,d\tau \leq c_j$ for all $s\in\R$,
\[
|y_1(\tau)-y_2(\tau)|\leq c_j\,e^{c_j}\,\sigma\quad \hbox{for almost every }\;\tau\in [a,b]\,.
\]
\end{lem}
\begin{proof}
Let $v(\tau) \coloneqq |y_1(\tau)-y_2(\tau)|=\left| \int_a^\tau g_1(y_1(s),s)\,ds -\int_a^\tau g_2(y_2(s),s)\,ds \right|$, $\tau\in [a,b]$. By the triangle inequality, using the $l$-bound $\wt \ell_j$ for $g_1$ and \eqref{eq:aux 4}, it is easy to obtain $v(\tau)\leq \int_a^\tau \wt\ell_j(s)\,v(s)\,ds + c_j\,\sigma$ for a.e.~$\tau\in [a,b]$. Then, the result follows by applying a version of Gronwall's lemma for $\wt \ell_j \in L^1_{loc}$ and the hypothesis of bounded $l$-bounds.
\end{proof}
The main result is the following.
\begin{thm}\label{teor:main}
Under Assumptions~{\rm\ref{asu:1}} and~{\rm\ref{asu:2}}, for each $\ep>0$  let $(x_\ep(t),y_\ep(t))$ be the solution of~\eqref{eq:slowfastent}-\eqref{eq:initial values}. Then, given $t_0>0$, for $\ep$ small enough, $(x_\ep(t),y_\ep(t))$ can be extended to $t\in [0,t_0]$. Furthermore:
\begin{itemize}
\item[\rm{(i)}]  Every sequence $\ep_j\to 0$ has a subsequence, say $\ep_k\downarrow 0$, such that $x_{\ep_k}(t)$ converges uniformly for $t\in [0,t_0]$ to a solution of the differential inclusion
\begin{equation*}
\dot{x}  \in \left\{ \int_{\R^m} f(x,y)\,d\mu(y) \left| \;\mu\in M\big(\A^x_{\Om(\bar\nu_g)}\big)\right.\right\},\quad x(0)=x_0\,,
\end{equation*}
where $M\big(\A^x_{\Om(\bar\nu_g)}\big)$ stands for the collection of projections on $\R^m$ of $\pi^x$-invariant measures supported on the attractor $\A^x$ restricted to the omega-limit set of $\bar\nu_g$, that is, $\A^x_{\Om(\bar\nu_g)}:=\A^x\cap (\Om(\bar\nu_g)\times \R^m)$. We denote by $x_0(t)$ the uniform limit of $x_{\ep_k}(t)$ for $t\in [0,t_0]$.
\item[\rm{(ii)}] For the fast variables in the fast timescale $y_{\ep_k}(\tau)$, the behaviour is as follows. For a fixed small $r>0$, let $K_0$ be the tubular set
\begin{equation*}
K_0:=\big\{x_0(t)+x\mid t\in [0,t_0],\;x\in\R^n \;{\rm with}\;|x|\leq r\big\}\subset \R^n
\end{equation*}
which is compact and let $\A\subset \mathcal{H} \times K_0\times \R^m$ be the global attractor given in Proposition~{\rm\ref{prop:existe atractor}} for the  skew-product semiflow~\eqref{eq:sk-pr}.  Then, given $\delta >0$, there exist $T=T(\delta,y_0,\A)>0$ and an integer $k_0=k_0(\delta,T)>0$ with $2T<t_0/\ep_{k_0}$ such that for every $k\geq k_0$,
\begin{equation}\label{eq:teorema}
 {\rm dist}\big(y_{\ep_k}(\tau),A_{(\bar\nu_g{\cdot}\tau,x_0(\ep_k\tau))}[\delta]\big)\leq  \delta \quad \hbox{for all }\; \tau\in [T,t_0/\ep_k]\,.
\end{equation}
In the special case where the global attractor $\A$ is a copy of the base, that is, when there is continuous map $\eta(\bar\nu,x)$ such that the sections $A_{(\bar\nu,x)}=\{\eta(\bar\nu,x)\}$ for all $(\bar\nu,x)\in\mathcal{H}\times K_0$, then, there is no need to consider the inflated version of the pullback attractors. That is: given $\delta >0$, there exist $T=T(\delta,y_0,\A)>0$ and an integer $k_0=k_0(\delta,T)>0$ with $2T<t_0/\ep_{k_0}$ such that for every $k\geq k_0$,
\begin{equation}\label{eq:teorema copy}
 \big|y_{\ep_k}(\tau)-\eta{(\bar\nu_g{\cdot}\tau,x_0(\ep_k\tau))}\big|\leq  \delta \quad \hbox{for all }\; \tau\in [T,t_0/\ep_k]\,.
\end{equation}

\end{itemize}
\end{thm}
\begin{proof}
The existence of the solutions $(x_\ep(t),y_\ep(t))$, $t\in [0,t_0]$ for $\ep$ small enough, and item (i) on the behaviour of the slow variables are part of Theorem~5.3 in~\cite{paper:ZA99}. The proof of (ii) essentially follows the steps of the proof of~\cite[Theorem~4.2]{paper:LOStracking}, so we simply indicate some small but necessary changes. Given a $\delta>0$, the way to determine the value of $T$ is the same, using Theorem~\ref{teor:forward} applied to a suitable compact subset $D\subset\R^m$ described in~\cite{paper:LOStracking}. Some differences arise when it comes to finding the value of the integer $k_0$. Associated with the compact sets $K_0\subset \R^n$ and $K:=y([0,2T]\times\mathcal{H}\times K_0\times D)\subset \R^m$, we find an integer $j\geq 1$ such that $K_0\times K\subset B_j^{(n)}\times B_j^{(m)}$. Next, we consider the $l$-bound $\wt l_j$ of $g$ and a constant $c_j>0$ such that $\int_s^{s+2T} \wt l_j(\tau)\,d\tau \leq c_j$ for all $s\in\R$, by the hypothesis of bounded $l$-bounds. Then, we take $\sigma>0$ small enough so that $\sigma\,c_j\,e^{c_j}\leq \delta/2$. For the modulus of continuity $\omega_j$ in Assumption~\ref{asu:1}~(ii)-(b), by the uniform convergence of $x_{\ep_k}$ to $x_0$ on $[0,t_0]$ and the fact that whenever $\tau,s\geq 0$ are such that $|\tau-s|\leq 2T$, $|\ep_k\tau-\ep_ks|\leq 2T\ep_k\downarrow 0$ as $k\to\infty$,  we can find a $k_0\geq 1$ such that $2T<t_0/\ep_{k_0}$ and for all $k\geq k_0$, $x_{\ep_k}(t)\in K_0$ for every $t\in [0,t_0]$ and $\omega_j(|x_{\ep_k}(\ep_k\tau) -x_{0}(\ep_ks)|)< \sigma$ for all $\tau,s\in [0,t_0/\ep_k]$ with $|\tau-s|\leq 2T$.

With this choice of $T$ and $k_0$, let us fix a $k\geq k_0$ and illustrate how to apply Lemma~\ref{lem:comparacion} in the first step of the finite iterative process to prove~\eqref{eq:teorema} on successive intervals $[T,2T]$, $[2T,3T],\ldots$ until the whole interval $[T,t_0/\ep_k]$ has been swept. The first step corresponds to the values of $\tau\in [T,2T]$. Fixed $s\in [T,2T]$,
\begin{itemize}
  \item $y_{\ep_k}(\tau)$ is a solution of $y'=g(x_{\ep_k}(\ep_k\tau),y,\tau)$ for $\tau\in [0,2T]$;
  \item the map $g_1(y,\tau):= g(x_{\ep_k}(\ep_k\tau),y,\tau)$ is a Carath\'{e}odory map  with $m$-bounds and bounded $l$-bounds, thanks to Assumption~\ref{asu:1};
  \item $y_2(\tau):=y(\tau,\bar\nu_g,x_0(\ep_k s),y_0)$ is a solution of $y'=g(x_0(\ep_k s),y,\tau)$ for $\tau\in [0,2T]$;
  \item $y_{\ep_k}(0)=y_2(0)=y_0\in D$, by the choice of $D$ (see~\cite{paper:LOStracking} for details);
  \item $y_{\ep_k}(\tau),\, y_2(\tau)\in K\subset B_j^{(m)}$  for every $\tau\in [0,2T]$, by construction;
  \item $\big|g(x_{\ep_k}(\ep_k\tau),y_2(\tau),\tau)-g(x_0(\ep_k s),y_2(\tau),\tau)\big|\leq \wt l_j(\tau)\,\omega_j(|x_{\ep_k}(\ep_k\tau) -x_{0}(\ep_ks)|) \\< \wt l_j(\tau)\,\sigma$ for a.e.~$\tau\in [0,2T]$, by Assumption~\ref{asu:1} and the choice of $k_0$.
\end{itemize}
Then, for $\tau\in [0,2T]$ Lemma~\ref{lem:comparacion} provides $|y_{\ep_k}(\tau)-y_2(\tau)|\leq \sigma\,c_j\,e^{c_j}\leq \delta/2$, by the choice of $\sigma$. Evaluating at $\tau=s$, this means that for every $s\in [T,2T]$,
\begin{equation}\label{eq:previa}
\big|y_{\ep_k}(s)-\left.y(\tau,\bar\nu_g,x_0(\ep_k s),y_0)\right|_{\tau=s}\big|\leq  \displaystyle\frac{\delta}{2}\,.
\end{equation}
Now, note that for $\tau\geq T$, since $y_0\in D$,  ${\rm dist}\big(y(\tau,\bar\nu_g,{x_0(\ep_k s)},y_0),A_{(\bar\nu_g{\cdot}\tau,x_0(\ep_k s))}[\delta]\big)$ $\leq {\rm dist}\big(y(\tau,\bar\nu_g,{x_0(\ep_k s)},D),A_{(\bar\nu_g{\cdot}\tau,x_0(\ep_k s))}[\delta]\big) <\delta/2$ by the choice of $T$. This, combined with~\eqref{eq:previa}, allows us to deduce that ${\rm dist}\big(y_{\ep_k}(s),A_{(\bar\nu_g{\cdot}s,x_0(\ep_k s))}[\delta]\big)\leq  \delta$ for every $s\in [T,2T]$, i.e.,~\eqref{eq:teorema} holds on the interval $[T,2T]$.
From this, the proof continues by adapting the arguments in the proof of~\cite[Theorem~4.2]{paper:LOStracking}, with no further difficulties. We just mention that the semicocycle property \eqref{eq:cocycle} is crucial in the following steps. The proof of \eqref{eq:teorema} is finished.

Improved versions of the last tracking result can be obtained under stronger assumptions such as the continuity of the setvalued maps $x\in K_0\mapsto \A^x\subset \mathcal{H}\times \R^m$ or  $(\bar\nu,x)\in  \mathcal{H}\times K_0\mapsto A_{(\bar\nu,x)}\subset \R^m$. The particular case where $\A$ is a copy of the base is specially useful in applications. We omit the proof of \eqref{eq:teorema copy}, since it follows similar arguments to the preceding ones. The reader is referred to Section~4 in~\cite{paper:LOStracking} for more details. The proof is finished.
\end{proof}
Finally, we stress that the extension of the results in~\cite{paper:LOStracking} to the weaker assumptions considered in this section is highly nontrivial. A key ingredient is the completion of the precompact family of translated ODEs, obtained by means of the theory of generalized ODEs associated with $N$-dimensional parametric b-measures developed in the preceding sections.

\subsection{Numerical examples.} To complete our study, we illustrate the tracking phenomenon using a scalar ordinary differential equation that explicitly depends on time and exhibits two separate time scales.
We consider the non-homogeneous linear problem,
\begin{equation}\label{eq:linear-cantor}
 y' = -\big(2 - \sin(\tau)\big)y + \cos(\ep \sqrt2 \tau)+ h(\tau)\,,
\end{equation}
where $h(\tau)$ is defined as in \eqref{eq:precompact map h}. Without the term $h(\tau)$ the problem is quasiperiodic, but its presence requires considering generalized equations in the construction of the hull.  It is immediate to see that this problem admits a unique hyperbolic solution. Indeed, given any two solutions $y(\tau,\tau_0,y_1)$ and $y(\tau,\tau_0,y_2)$, we have that $y(\tau,\tau_0,y_1)- y(\tau,\tau_0,y_2)$ satisfies the homogeneous equation
\[
 u' = -\big(2 - \sin(\tau)\big)u\,,
\]
whose general solution $u(\tau,\tau_0,u_0)$ satisfies
\[
|u(\tau,\tau_0,u_0)|=|u_0\,e^{\cos(\tau_0)-\cos(\tau)}\,e^{-2(\tau-\tau_0)}|\le |u_0|\,e^2\,e^{-2(\tau-\tau_0)}\xrightarrow{\tau\to\infty}0\,,
\]
which implies that any pair of trajectories converges exponentially to each other in forward time. On the other hand, note that \eqref{eq:linear-cantor} can be solved explicitly using the variation of constants formula. The solution with initial datum $y_0$ at time $\tau_0$ is
\begin{equation}\label{eq:linear-sol}
y(\tau,\tau_0,y_0)=a(\tau,\tau_0)\,y_0+\int_{\tau_0}^\tau a(\tau,s)\big(\cos(\ep \sqrt2 s)+ h(s)\big)\,ds
\end{equation}
where $a(\tau,s)=e^{\cos(s)-\cos(\tau)}e^{-2(\tau-s)}$. Hence, fixing $\tau\in\R$ and taking the pullback limit as $\tau_0\to-\infty$ in \eqref{eq:linear-sol}, one has the unique bounded solution:
\[
\overline y(\tau)= \int_{-\infty}^\tau a(\tau,s)\big(\cos(\ep \sqrt2 s)+ h(s) \big)\,ds\,,
\]
which is therefore the sought hyperbolic attracting solution. Following the same arguments, it is also immediate to show that each equation of the family
\begin{equation}\label{eq:linear-cantor-parametric}
y' = -\big(2 - \sin(\tau)\big)y + c+ h(\tau)\,,\qquad c\in[-1,1]\,,
\end{equation}
also has a unique hyperbolic attracting solution
\[
\overline y_c(\tau)= \int_{-\infty}^\tau a(\tau,s)\big(c+ h(s)\big)\,ds\,,\qquad c\in[-1,1]\,.
\]
Note also that \eqref{eq:linear-cantor} can be written in the form of a slow-fast system,
\[
\begin{cases}
     x'=\ep\\
     y' = -\big(2 + \sin(\tau)\big)y + \cos(\sqrt2 x)+ h(\tau)\,.
\end{cases}
\]
We remark that, due to the presence of the irregular forcing $h(\tau)$, the results in \cite{paper:LOStracking} do not apply. The reader will be able to check that the assumptions of Theorem~\ref{teor:main} are satisfied and the considered hyperbolic solutions for each parametric problem and for the slow-fast problem lift to the global attractors of the corresponding skew-product flows. In fact, thanks to the non-homogeneous linear structure, it is easy to check that all the equations in the hull also have a unique hyperbolic attracting solution, so that the global attractor $\A$ of the induced skew-product semiflow is a copy of the base. Then,  the tracking formula \eqref{eq:teorema copy} holds.
\begin{figure}
    \centering
    \begin{overpic}[width=\textwidth]{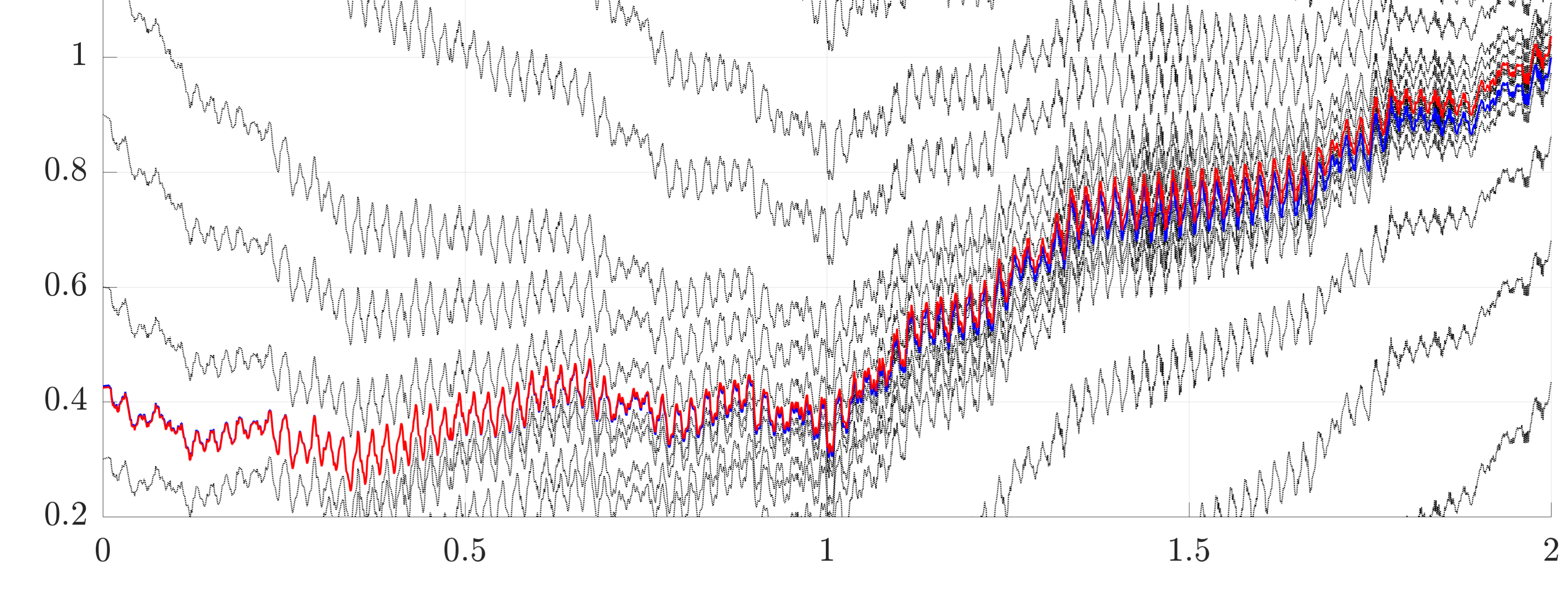}
    \put(50,-2){$\tau$}
    \put(7,35){$y$}
    \end{overpic}
    \caption{$\ep=0.15$. In red the hyperbolic attracting solution for \eqref{eq:linear-cantor}. In blue, the curve obtained by considering the value at time $\tau$ of the attractor of \eqref{eq:linear-cantor-parametric} obtained for $c=\cos(\ep\sqrt2\tau)$. Dotted black lines are sample solutions of \eqref{eq:linear-cantor} for selected initial conditions.}
    \label{fig:Cantor-tracking-1}
\end{figure}
In Figure \ref{fig:Cantor-tracking-1} it is possible to appreciate how solutions of the slow-fast equation approach the parametric graph in solid blue obtained depicting the value at time $\tau$ of the globally attracting hyperbolic solution to the parametric equation \eqref{eq:linear-cantor-parametric} corresponding to the value of the parameter $c=\cos(\varepsilon\sqrt{2}\tau)$. This is also true for the globally attracting hyperbolic solution of the slow-fast problem depicted in solid red on the interval $[0,1/\varepsilon]$ (depicted only in $[0,2]$ to maximize the image's clarity).\smallskip
\begin{figure}
    \centering
\begin{overpic}[abs,unit=0.5mm,scale=.25, width=0.48\textwidth]{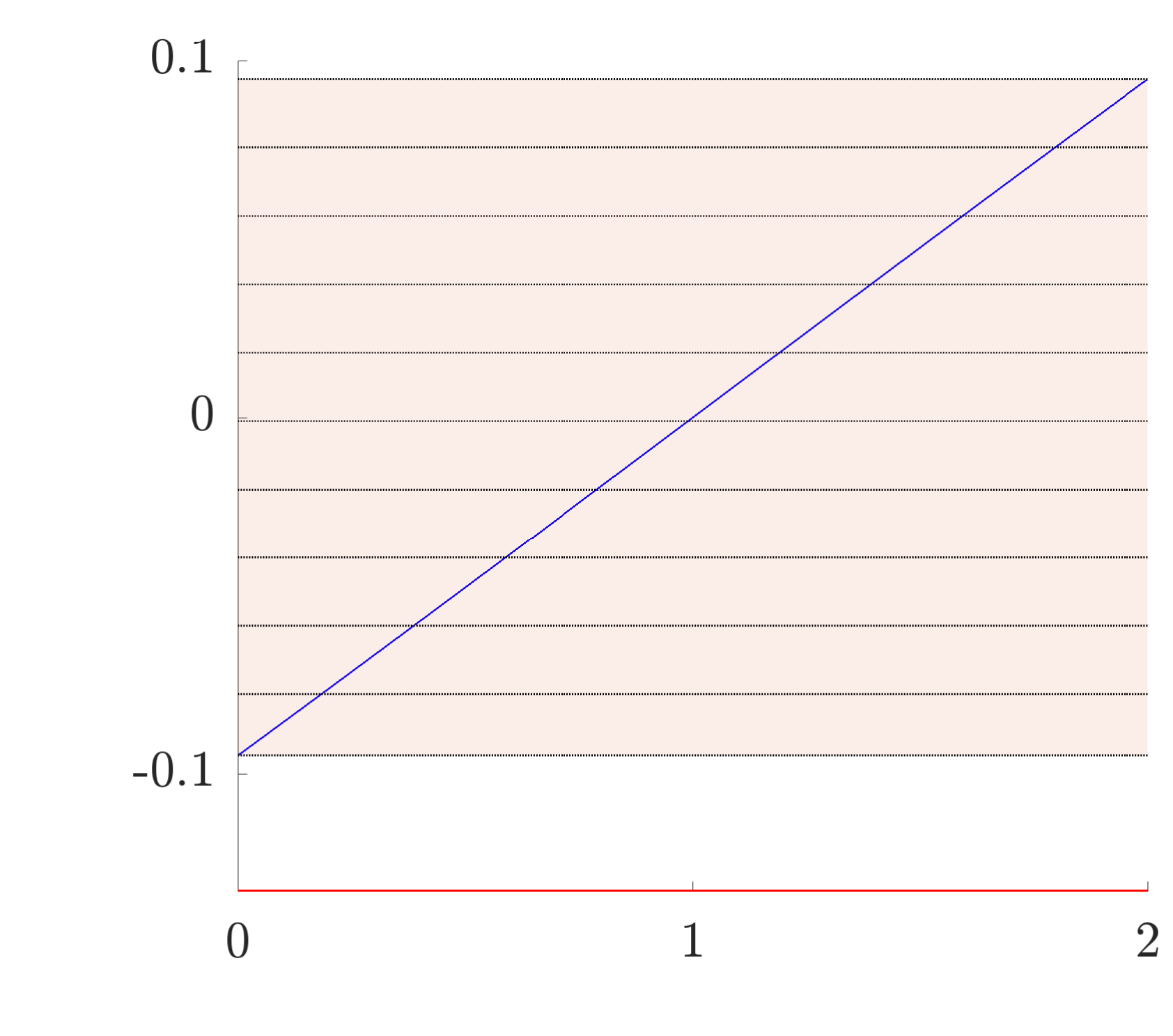}
\put(26,94){$c$}
\put(110,17){$\tau$}
\end{overpic}
\hspace{2mm}
\begin{overpic}[abs,unit=0.5mm,scale=.25,  width=0.48\textwidth]{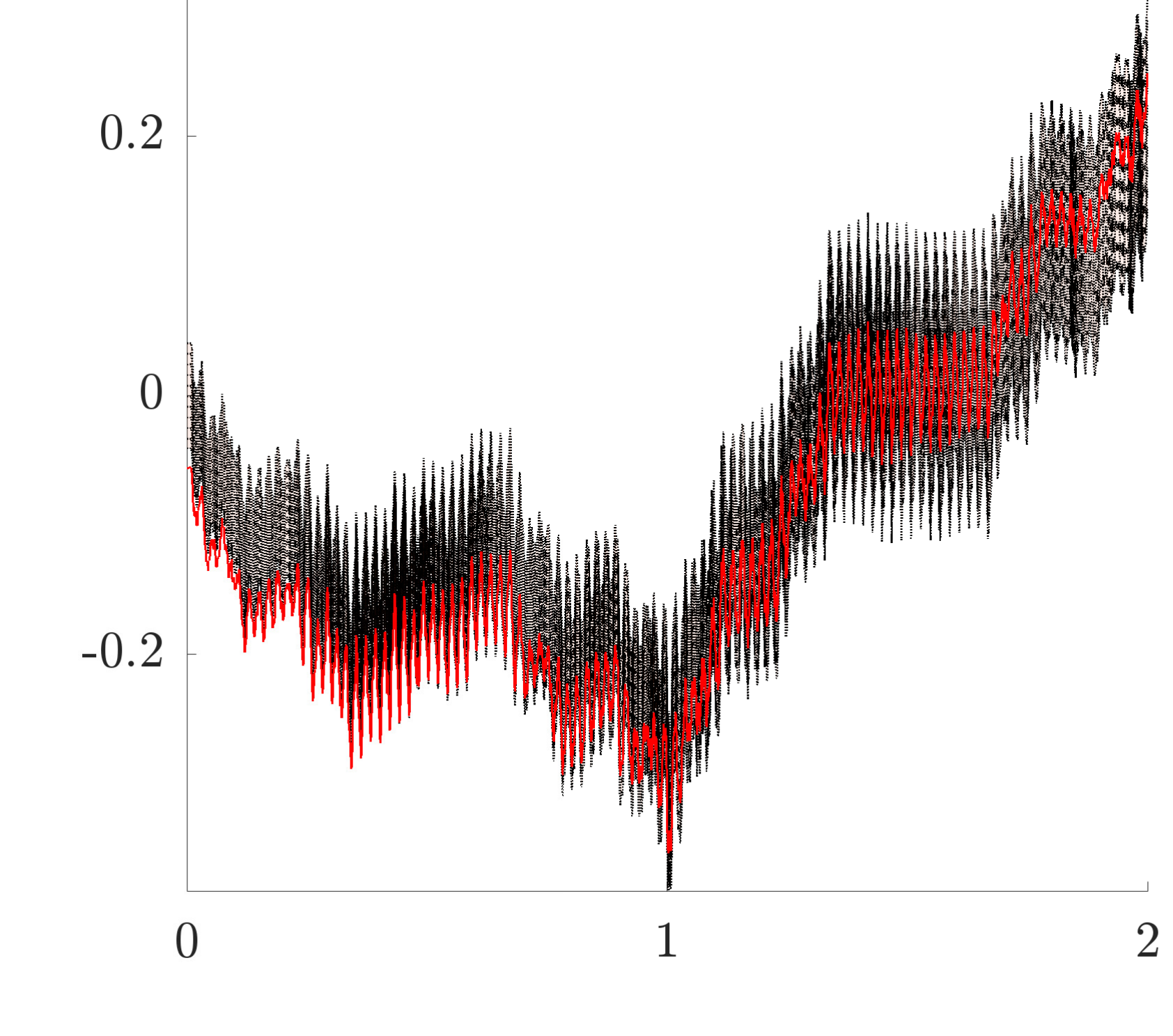}
\put(22,94){$y$}
\put(110,17){$\tau$}
\end{overpic}\\
\vspace{5mm}
\begin{overpic}[trim={0 1cm 0 2cm},clip,  width=\textwidth]{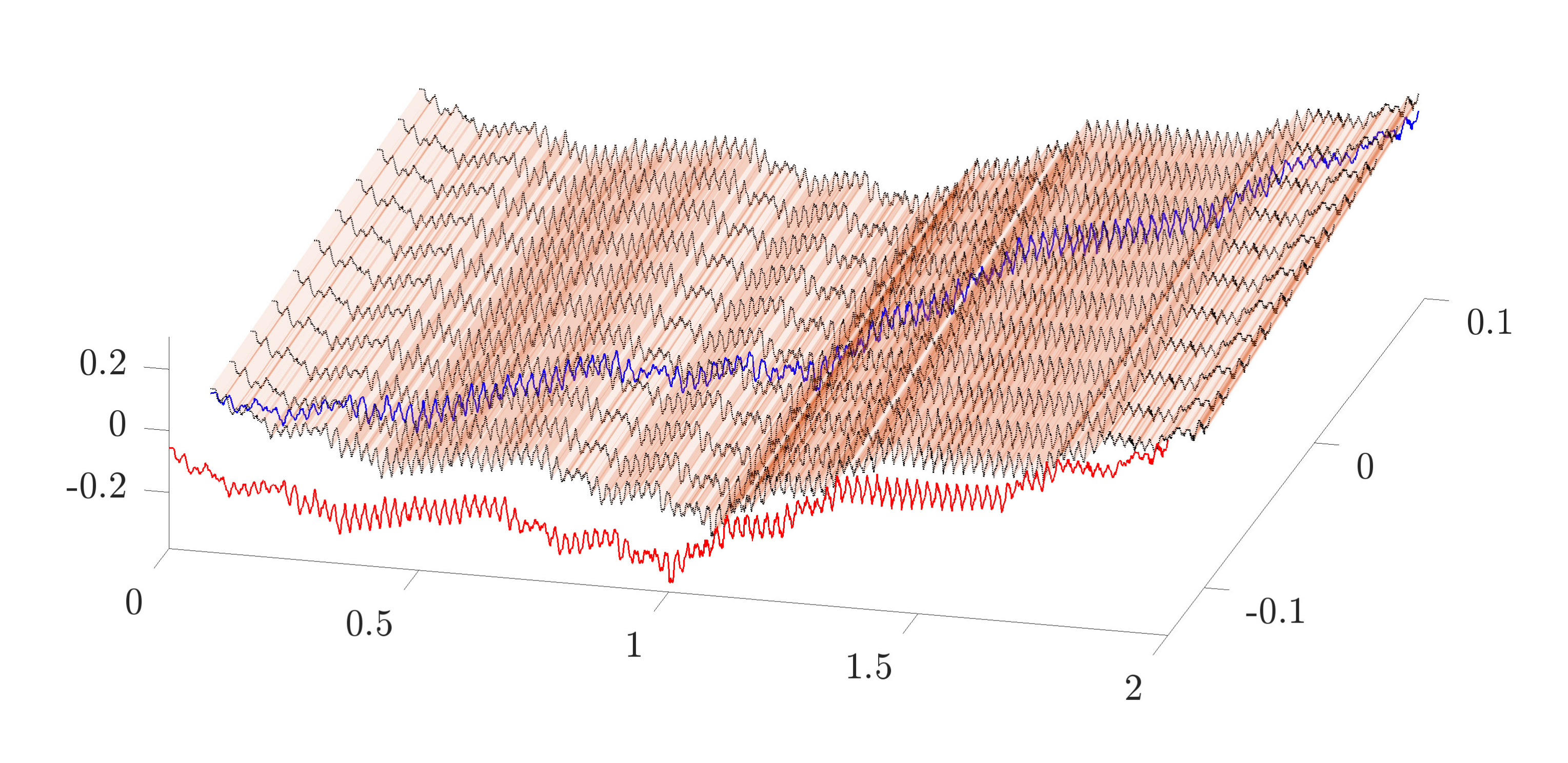}
\put(2,19){$y$}
\put(39,2){$\tau$}
\put(90,14){$c$}
\end{overpic}

\hspace{2mm}

\caption{$\ep=0.15$. Numerical representation of the implications of Theorem \ref{teor:main} in the special case where the global attractor is a copy of the base. The three panels correspond to the same plot from three different points of view. The orange surface is the continuum of copies of the base of \eqref{eq:linear-cantor-parametric} for the values of $c$ in the image of the function $2\arctan(s)/\pi$; the blue line is the graph of the function $\eta{(\bar\nu_g{\cdot}\tau,\ep\tau)}$, made taking at time $\tau$ the value of the copy of the base at time $\tau$ obtained for $c=2\arctan(\varepsilon \tau)/\pi$---it is possible to appreciate how the blue curve follows the variation of the attractors of the parametric family \eqref{eq:linear-cantor-parametric} on the slow time-scale; the red line is the attracting hyperbolic solution of \eqref{eq:linear-cantor-2}. The top right panel---projection on $(\tau,y)$ shows how the hyperbolic solution of \eqref{eq:linear-cantor-2} perfectly overlaps $\eta{(\bar\nu_g{\cdot}\tau,\ep\tau)}$ on $[0,2]$.}
\label{fig:Cantor-tracking-2}
\end{figure}

To further emphasize the relationship between the globally attracting solutions of the  parametric problem \eqref{eq:linear-cantor-parametric} and the one of the slow-fast equation we consider a modification of \eqref{eq:linear-cantor} in that we use the function $2\arctan(\ep \tau)/\pi$ instead of $\cos(\ep\sqrt 2 \tau)$, which yields the differential equation
\begin{equation}\label{eq:linear-cantor-2}
 y' = -\big(2 - \sin(\tau)\big)y + \frac{2}{\pi}\arctan(\ep \tau)+ h(\tau)\,.
\end{equation}
The graph of the arc-tangent is almost linear in the interval $[0,2]$ (cf.~the first panel of Figure \ref{fig:Cantor-tracking-2}) and the hyperbolic solutions of the equations \eqref{eq:linear-cantor-parametric} for $c$ fixed in the image of $2\arctan(s)/\pi$ are the same up to an upward translation. They define a continuum of hyperbolic solutions with variation on the fast time-scale (cf.~orange surface in the lower panel of Figure \ref{fig:Cantor-tracking-2}). It is therefore possible to draw the graph of the function $\eta{(\bar\nu_g{\cdot}\tau,\ep\tau)}$ which is made taking at time $\tau$ the value of the copy of the base at time $\tau$ obtained for $c=2\arctan(\varepsilon \tau)/\pi$ (blue curve in the lower panel of Figure \ref{fig:Cantor-tracking-2}). The reader will be able to appreciate that on the slow time-scale, this curve is actually crossing the continuum of hyperbolic solutions of the parameteric family \eqref{eq:linear-cantor-parametric} (orange surface). The red solid line is instead the hyperbolic solution of equation \eqref{eq:linear-cantor-2}. The right-hand top panel in Figure \ref{fig:Cantor-tracking-2} shows that the red line indistinguishably overlaps the blue one.\smallskip

All the above examples are simulated in Matlab R2023b using the closed formula for the hyperbolic solution of a scalar linear problem (as explicitly obtained above). The integral is approximated with the trapezoidal method {\rm trapz} starting at $t_0=-20$. 

\end{document}